\documentclass[11pt]{article}
\usepackage{amsfonts}
\usepackage{mathrsfs,color,amscd,amssymb,enumerate,amsthm,amsmath,bm,graphicx,psfrag,subfigure}
\usepackage{latexsym}
\usepackage{float,fancybox,shapepar,setspace,hyperref}
\usepackage{pgf,tikz}

\makeatletter
\def\leftharpoonfill@{\arrowfill@\leftharpoonup\relbar\relbar}
\def\rightharpoonfill@{\arrowfill@\relbar\relbar\rightharpoonup}
\newcommand\rbjt{\mathpalette{\overarrow@\rightharpoonfill@}}
\newcommand\lbjt{\mathpalette{\overarrow@\leftharpoonfill@}}
\makeatother

\makeatletter

\renewcommand{\@seccntformat}[1]{{\csname the#1\endcsname}{\normalsize .}\hspace{.5em}}
\makeatother

\def \[{\begin{equation}}
\def \]{\end{equation}}

\def \ss {\subseteq}

\def\bqed{ \hfill $\blacksquare$}

\newtheorem{Question}{Question}

\newtheorem{thm}{Theorem}

\newtheorem{claim}{Claim}

\newtheorem{lem}{Lemma}

\usetikzlibrary{arrows}
\begin{document}

\title{Distribution of  cycle lengths in graphs}
\author{Xiaolin Wang\footnote{Corresponding author. Email: xiaolinw@fzu.edu.cn}, Meiduo Chen, Xueping Xu\\
\small{School  of  Mathematics  and  Statistics,  Fuzhou  University,  Fuzhou  350108,  P.R. CHINA}\\
}
\date{}
\maketitle
\begin{abstract}

We present the relations between clique number and chromatic number with given the number of odd or even or all cycle lengths. 
Let $L_o(G)$ be the set of odd cycle lengths of $G$ and  $\ell_o(G)$ be the longest odd cycle length.  Gy\'arf\'as (DM, 1992) proved a classic result: $\chi(G)\leq 2|L_o(G)|+2$, and if $w(G)\leq 2|L_o(G)|+1$, then $\chi(G)\leq 2|L_o(G)|+1$. Later, Wang (SIAM DM, 2008), Ma and Ning (SIAM DM, 2018) together determined the exact chromatic number when $|L_o(G)|=2$. 
We further prove that if 
$w(G)\leq 2|L_o(G)|$ for any  $|L_o(G)|\geq 2$, then $\chi(G)\leq 2|L_o(G)|$. We also construct a class of graphs with $w(G)=2|L_o(G)|-1$ but $\chi(G)=2|L_o(G)|$ for every $|L_o(G)|\geq 2$. Using our result, we give a short proof of the relation between $w(G)$ and $\chi(G)$ with given $\ell_o(G)$ proved by Kenkre and  Vishwanathan (JGT, 2006).

Let $L_e(G)$ be the set of even cycle lengths of $G$ and  $\ell_e(G)$ be the longest even cycle length. Mihók and  Schiermeyer (DM, 2004) proved that $\chi(G)\leq 2|L_e(G)|+3$, and if $w(G)\leq 2|L_e(G)|+2$, then $\chi(G)\leq 2|L_e(G)|+2$. 
We further prove that if 
$w(G)\leq 2|L_e(G)|+1$ and  $|L_e(G)|\geq 3$, then $\chi(G)\leq 2|L_e(G)|+1$. We also construct a class of graphs with $w(G)=2|L_e(G)|$ but $\chi(G)=2|L_e(G)|+1$ for every $|L_e(G)|\geq 2$. Our result can deduce  the relation between $w(G)$ and $\chi(G)$ with given $\ell_e(G)$. 

Combining all the above results, we deduce the similar relations between  $w(G)$ and $\chi(G)$ with given the number of cycle lengths or the longest cycle length.
\vskip 2mm

\end{abstract}
{\setcounter{section}{0}

\section{Introduction}\setcounter{equation}{0}
\vskip 2mm
Let $G=(V(G),E(G))$ be a  simple undirected connected graph. Let $L_o(G)$ denote the set of odd cycle lengths of $G$, that is, $$L_o(G)=\{2i+1~|~G~contains~an~odd~cycle~of~length~2i+1\}.$$
If $G$ has no odd cycle, then  $|L_o(G)|=0$.
Let $\ell_o (G)$ be the length of the longest odd cycle in a non-bipartite graph $G$. Obviously, if $G$ is non-bipartite, then $2|L_o(G)|+1\leq \ell_o(G)$. And the equality holds if and only if $G$ contains  all odd cycles of  lengths  at most $ \ell_o(G)$. 
Similarly, let $L_e(G)$ denote the set of even cycle lengths of $G$, that is, $$L_e(G)=\{2i+2~|~G~contains~an~even~cycle~of~length~2i+2\}.$$
If $G$ has no even cycle, then  $|L_e(G)|=0$.
Let $\ell_e (G)$ be the length of the longest even cycle in  $G$. Obviously, if $|L_e(G)|\geq 1$, then $2|L_e(G)|+2\leq \ell_e(G)$. And the equality holds if and only if $G$ contains  all even cycles of  lengths  at most $ \ell_e(G)$. Let $L(G)$ denote the set of all cycle lengths of $G$ and  $\ell(G)$ denote the longest cycle length of $G$. Then $|L(G)|+2\leq \ell(G)$. 
Denote the clique number $w(G)$ as  the size of the maximum clique in $G$. Let $\delta(G)$ be the minimum degree in $G$.

In this paper, we study the relation between cycle lengths, $ \delta(G),$ $w(G)$ and $\chi(G)$,  which is a fundamental and extensive research area in graph theory.
In 1990, Bollob\'as and Erd\H{o}s \cite{E} conjectured that $\chi(G)\leq 2|L_o(G)|+2$ for any graph $G$. Two years later,
 Gy\'arf\'as 
confirmed this conjecture   and proved a strong result. Note that $w(G)\leq 2|L_o(G)|+2$.

\begin{thm}(Gy\'arf\'as \cite{G})\label{AG}
    $\chi(G)\leq 2|L_o(G)|+2$. Moreover, if $w(G)\leq 2|L_o(G)|+1$, then $\chi(G)\leq 2|L_o(G)|+1$.
\end{thm}

By Theorem \ref{AG},  $\chi(G)= 2|L_o(G)|+2$ if and only if $w(G)=2|L_o(G)|+2$. 
If one considers the elements of $L_o(G)$, then Gy\'arf\'as' result can be improved. If $|L_o(G)|=0$, then $G$ is bipartite. If $|L_o(G)|=1$, then by Theorem \ref{AG}, $\chi(G)=\max\{3,w(G)\}$. 
If $|L_o(G)|=2$, the problem becomes  hard.
Define a \emph{wheel} $W_n$ as a graph obtained by
adding all edges between $C_n$ and a new vertex.  The first result to $|L_o(G)|=2$ was from Wang \cite{W}.

\begin{thm}(Wang \cite{W})\label{Wang}
    Let $L_o(G)=\{3,5\}$. If $G$ contains neither $K_4$ nor $W_5$, then $\chi(G)=3$. Otherwise, $\chi(G)=\max\{4,w(G)\}$.
\end{thm}

Not long after, Camacho and  Schiermeyer \cite{CS} improved the bound to $\chi(G)\leq 4$ when $L_o(G)=\{3+2s,5+2s\}$ for any positive integer $s$.
The case $L_o(G)=\{5,7\}$ was determined to  $\chi(G)=3$ by
Kaiser,  Ruck\'y and R. $\check{\text{S}}$krekovski \cite{KR}. A significant breakthrough was made in  2018. Ma and Ning \cite{MN} proved the following theorem.

\begin{thm}(Ma and Ning \cite{MN})\label{Majie}
  Let $s,t$ be two positive integer.   If $L_o(G)=\{3,5+2s\}$, then $\chi(G)=\max\{3,w(G)\}$.  If $L_o(G)=\{3+2t,3+2t+2s\}$, then $\chi(G)=3$.
\end{thm}

Together with the above results, these researchers gave a complete solution to $|L_o(G)|=2$. Although it is very difficult to determine the chromatic number when $|L_o(G)|\geq 3$, it is also interesting to improve Gy\'arf\'as' result with smaller clique number. 
In this paper, 
we first prove the following theorem.

\begin{thm}\label{2connected}
   Let $G$ be a 2-connected non-bipartite graph with $\delta(G)\geq 2k$ and $k\geq 3$. Then $|L_o(G)|\geq k$. Moreover, if $|L_o(G)|=k$, then $2|L_o(G)|+1=\ell_o(G)$, and either $K_{2k+1}\ss G$ or $\chi(G)\leq 2k$.
\end{thm}

By Theorems \ref{AG} - \ref{2connected}, we can deduce the following theorem.

\begin{thm}\label{nok2L+1}
   Let $G$ be a graph with $|L_o(G)|\geq 2$. If  $w(G)\leq 2|L_o(G)|$, then  $\chi(G)\leq 2|L_o(G)|$. 
\end{thm}

By Theorems \ref{AG} and \ref{nok2L+1},  when $|L_o(G)|\geq 2$, $\chi(G)= 2|L_o(G)|+1$ if and only if $w(G)=2|L_o(G)|+1$. For the longest odd cycle length $\ell_o(G)$, we have some similar results. By Theorem \ref{AG},  $\chi (G)\leq \ell_o(G)+1$, which was also proved by Erd\H{o}s and Hajnal \cite{We} independently.
Applying Theorems \ref{AG} and  \ref{nok2L+1}, we give a short proof of two results of Kenkre and  Vishwanathan \cite{KV}.

\begin{thm}(Kenkre and  Vishwanathan \cite{KV})\label{KAV}
 Let $G$ be a  graph with $|L_o(G)|\geq 2$.   If $w(G)\leq \ell_o(G)$, then $\chi (G)\leq \ell_o(G)$. If $w(G)\leq \ell_o(G)-1$, then $\chi (G) \leq \ell_o(G)-1$.
\end{thm}

\noindent{\bf{Proof.}}
If $\ell_o(G)=2|L_o(G)|+1$, then these two results hold by Theorems \ref{AG} and \ref{nok2L+1}. Otherwise, $\ell_o(G)\geq 2|L_o(G)|+3$ and then $w(G)\leq 2|L_o(G)|$. By Theorem \ref{nok2L+1}, we get $\chi(G)\leq 2|L_o(G)|\leq \ell_o(G)-3$. \bqed

\vskip 2mm
By Theorem \ref{KAV}, when $|L_o(G)|\geq 2$,
$\chi(G)= \ell_o(G)+i$ if and only if $w(G)=\ell_o(G)+i$ for any $i\in \{0,1\}$.
By the proof of Theorem \ref{KAV}, we can deduce that $\ell_o(G)=2|L_o(G)|+1$ is a necessary condition for $\chi(G)=\ell_o(G)-1$.

One may ask if $w(G)\leq 2|L_o(G)|-1$ or if $w(G)\leq \ell_o(G)-2$, can we have $\chi(G)\leq 2|L_o(G)|-1$ or $\chi(G)\leq \ell_o(G)-2$? They are both not true by the following example. Let   $H_i$ be a graph obtained by adding all edges between  $K_{i}$ and $C_5$. It is easy to see that if $i=2k-3\geq 1$, then $w(H_{2k-3})=2k-1$ and $2|L_o(H_{2k-3})|+1=\ell_o(H_{2k-3})=2k+1$. But we can check that  $\chi(H_{2k-3})=2k=2|L_o(H_{2k-3})|=\ell_o(H_{2k-3})-1$. Since this counterexample  satisfies $2|L_o(G)|+1=\ell_o(G)$, it is interesting to consider the following question.

\begin{Question}
Let $G$ be a graph with $\ell_o(G)>2|L_o(G)|+1$.
Does $\chi(G)\leq 2|L_o(G)|-1$ follow from 
$w(G)\leq 2|L_o(G)|-1$, or does $\chi(G)\leq \ell_o(G)-2$ follow from $w(G)\leq \ell_o(G)-2$?

\end{Question}

The second part of this question is true because $\ell_o(G)=2|L_o(G)|+1$ is a necessary condition for $\chi(G)=\ell_o(G)-1$.
When $|L_o(G)|=2$, by Theorem \ref{Majie}, the first part of this question  is true. So we leave the first part of this question  for further research.
It is natural to ask, when considering even cycle lengths, can we have  analogous results?
In 2004, Mihók and  Schiermeyer first \cite{MS} proved the analogous result for even cycle lengths. Note that $w(G)\leq 2|L_e(G)|+3$.

\begin{thm}(Mihók and  Schiermeyer \cite{MS})\label{MS}
    $\chi(G)\leq 2|L_e(G)|+3$. Moreover, if $w(G)\leq 2|L_e(G)|+2$ and $|L_e(G)|\geq 1$, then $\chi(G)\leq 2|L_e(G)|+2$.
\end{thm}

When $|L_e(G)|=0$, by Theorem \ref{MS}, $\chi(G)\leq 3$. If $|L_e(G)|\geq 1$, by Theorem \ref{MS},  $\chi(G)=2|L_e(G)|+3$ if and only if $w(G)=2|L_e(G)|+3$.
In this paper, we prove the following  theorem similar to Theorem \ref{2connected}.

\begin{thm}\label{2connectedeven}
   Let $G$ be a 2-connected  graph with $\delta(G)\geq 2k+1$ and $k\geq 3$. Then $|L_e(G)|\geq k$. Moreover, if $|L_e(G)|=k$, then $2|L_e(G)|+2=\ell_e(G)$, and either $K_{2k+2}\ss G$ or $\chi(G)\leq 2k+1$. 
\end{thm}

Applying Theorems \ref{MS} and \ref{2connectedeven}, we have the following theorem.

\begin{thm}\label{nok2L+2}
   Let $G$ be a graph with $|L_e(G)|\geq 3$ and $w(G)\leq 2|L_e(G)|+1$. Then we have $\chi(G)\leq 2|L_e(G)|+1$. 
\end{thm}

By Theorems \ref{MS} and \ref{nok2L+2},  $\chi(G)=2|L_e(G)|+2$ if and only if $w(G)=2|L_e(G)|+2$ when $|L_e(G)|\geq 3$.
Applying   Theorems \ref{MS} and  \ref{nok2L+2}, we can also  have the following result for the longest even cycle length.

\begin{thm}\label{longesteven}
$\chi(G)\leq \ell_e(G)+1$.  If $w(G)\leq \ell_e(G)$ and $|L_e(G)|\geq 3$, then $\chi (G)\leq \ell_e(G)$. If $w(G)\leq \ell_e(G)-1$ and $|L_e(G)|\geq 3$, then $\chi (G) \leq \ell_e(G)-1$.
\end{thm}

\noindent{\bf{Proof.}}
If $\ell_e(G)=2|L_e(G)|+2$, then the  results hold by Theorems \ref{MS} and \ref{nok2L+2}. Otherwise, $\ell_e(G)\geq 2|L_e(G)|+4$ and then $w(G)\leq 2|L_e(G)|+1$. By Theorem \ref{nok2L+2}, we get $\chi(G)\leq 2|L_e(G)|+1\leq \ell_e(G)-3$. \bqed
\vskip 2mm

By Theorem \ref{longesteven}, when $|L_e(G)|\geq 3$,
$\chi(G)= \ell_e(G)+i$ if and only if $w(G)=\ell_e(G)+i$ for any $i\in \{0,1\}$.
By the proof of Theorem \ref{longesteven}, we can deduce that $\ell_e(G)=2|L_o(G)|+2$ is a necessary condition for $\chi(G)=\ell_e(G)-1$.

One may ask if $w(G)\leq 2|L_e(G)|$ or if $w(G)\leq \ell_e(G)-2$, can we have $\chi(G)\leq 2|L_e(G)|$ or $\chi(G)\leq \ell_e(G)-2$? 
It is easy to see that if $2k-2\geq 2$, then $w(H_{2k-2})=2k$ and $2|L_e(H_{2k-2})|+2=\ell_e(H_{2k-2})=2k+2$. But we can check that  $\chi(H_{2k-2})=2k+1=2|L_e(H_{2k-2})|+1=\ell_e(H_{2k-2})-1$.
Since this counterexample  satisfies $2|L_e(G)|+2=\ell_e(G)$, it is interesting to consider the following question.

\begin{Question}
Let $G$ be a graph with $\ell_e(G)>2|L_e(G)|+2$. Does $\chi(G)\leq 2|L_e(G)|$ follow from $w(G)\leq 2|L_e(G)|$, or does $\chi(G)\leq \ell_e(G)-2$ follow from $w(G)\leq \ell_e(G)-2$?
\end{Question}

The second part of this question is true because $\ell_e(G)=2|L_e(G)|+2$ is a necessary condition for $\chi(G)=\ell_e(G)-1$.
 So we leave the first part of this question  for further research.

One may concerning the chromatic number of $G$ with given $L(G)$. If $|L(G)|\leq 3$, we can obtain the exact value of 
$\chi(G)$ by considering $|L_o(G)|$ and $|L_e(G)|$, and by Theorems \ref{AG}, \ref{Wang}, \ref{Majie} and \ref{MS}. It is also interesting to consider the
relation between $w(G)$ and $\chi(G)$ with given $|L(G)|\geq 4$.  Note that $w(G)\leq |L(G)|+2$.

\begin{thm}\label{LG}
  Let $|L(G)|\geq 4$.  For each $i\in \{0,1,2\}$, if $w(G)\leq |L(G)|+i$, then $\chi(G)\leq |L(G)|+i$, except when $i=0$, $|L_o(G)|=3$ and $|L_e(G)|=2$.
    \end{thm}

\noindent{\bf{Proof.}} We first suppose that $|L(G)|=2k$. Then  $k\geq 2$. If $|L_o(G)|\leq k$, then the conclusion holds by Theorems \ref{AG} and \ref{nok2L+1}. Otherwise,  
$|L_e(G)|\leq k-1$. By Theorem \ref{MS}, we can check that the conclusion holds.

Suppose $|L(G)|=2k+1$.  Then  $k\geq 2$. If $|L_o(G)|\leq k$, then the conclusion holds by Theorem \ref{AG}. Otherwise,  
$|L_e(G)|\leq k$. Note that when $i=0$ and $|L(G)|=2k+1=5$, $|L_e(G)|\not=2$.   By Theorems \ref{MS} and \ref{nok2L+2}, we can check that the conclusion holds. \bqed

\vskip 2mm

 Note that $w(G)\leq |L(G)|+2\leq \ell(G)$.
We have the following theorem for $\ell(G)$.

\begin{thm}\label{lG}
    Let $|L(G)|\geq 4$.  For each $i\in \{0,1,2\}$, if $w(G)\leq \ell(G)-i$, then $\chi(G)\leq \ell(G)-i$, except when $i=2$, $|L_o(G)|=3$ and $|L_e(G)|=2$.
\end{thm}

\noindent{\bf{Proof.}} 
If $|L(G)|+2= \ell(G)$, then the conclusion holds by Theorem \ref{LG}. 
If 
$|L(G)|+3\leq  \ell(G)$, then $w(G)\leq |L(G)|+1$. By Theorem \ref{LG}, 
$\chi(G)\leq |L(G)|+1\leq \ell(G)-2$. 
\bqed
\vskip 2mm

We believe that the remaining case in Theorems \ref{LG} and \ref{lG} is true. And we leave this problem for further research. 
For $k\geq 4$, $|L(H_{k-3})|+2=\ell(H_{k-3})=k+2$, $w(H_{k-3})=k-1$ and $\chi(H_{k-3})=k$. So $\chi(G)\leq |L(G)|-1$  cannot follow from $w(G)\leq |L(G)|-1$, and $\chi(G)\leq \ell(G)-3$  cannot follow from $w(G)\leq \ell(G)-3$. 
Since this counterexample  satisfies $|L(G)|+2=\ell(G)$, it is interesting to consider the following question.

\begin{Question}
Let $G$ be a graph with $\ell(G)>|L(G)|+2$. Does $\chi(G)\leq |L(G)|-1$ follow from $w(G)\leq |L(G)|-1$, or does $\chi(G)\leq \ell(G)-3$ follow from $w(G)\leq \ell(G)-3$?
\end{Question}


In Theorems \ref{2connected} or \ref{2connectedeven}, we characterize the graphs with exactly $k$ odd or even cycle lengths, respectively. So, Theorems \ref{2connected} and \ref{2connectedeven} also improve the following two results.

\begin{thm}\label{LWZ}(Lin, Wang and Zhou \cite{LWZ})
If $G$ is a 2-connected non-bipartite graph with $\delta(G)\geq 2k$, then $G$ contains $k$ cycles with consecutive odd lengths. 
\end{thm}

\begin{thm}
    \label{LM}(Liu and Ma \cite{LM})
If $G$ is a 2-connected  graph with $\delta(G)\geq 2k+1$, then $G$ contains $k$ cycles with consecutive even lengths.
\end{thm}

These two theorems will help us to prove Theorems \ref{2connected} and \ref{2connectedeven} (although we can prove them without using these two theorems).
For more results about the relation between minimum degree, chromatic number  and cycle lengths, one can see  \cite{Ch}, \cite{Co}, \cite{Di}, \cite{GHL}, \cite{GHM}, \cite{GL}, \cite{GKS}, \cite{KSV},  \cite{LR}, 
\cite{Lo},
\cite{T}. For directed version of this topic, see \cite{B} and \cite{Ch}.

In the end of this section, we define some notations we will use later.
Denote the number of edges of a path $P$, i.e., the length of  $P$ as $|P|$. Let $C$ be a cycle. In this paper, we will define one direction of $C$ as the forward direction  and denote as $\overrightarrow{C}$.  And another as the backward direction  and denote as $\overleftarrow{C}$.  For any $x\in C$, denote $x^+,x^-\in C$ as the successor and predecessor of $x\in C$ in the forward direction, respectively. A cycle is called a $k$-cycle if $|C|=k$. For convenience, we simply denote $a\equiv b \pmod{2}$ or $a\not\equiv b \pmod{2}$ as $a\equiv b$ or $a\not\equiv b$, respectively. Let $K_n^{2-}$ ($n\geq 4$) be a graph obtained by deleting two disjoint edges in $K_n$. For other notations not introduced here, we refer to \cite{BM}.

The remain of this paper is organized as follows.   
In Section 2, we use Theorems \ref{2connected} and \ref{2connectedeven} to prove Theorems \ref{nok2L+1} and \ref{nok2L+2}, respectively. We will prove Theorem \ref{2connected} in Section 3 and Theorem \ref{2connectedeven} in Section 4. For easy to follow, we prove some important lemmas in Section 5.

\section{Proof of Theorems \ref{nok2L+1} and \ref{nok2L+2}}

\noindent{\bf{Proof of Theorem \ref{nok2L+1}}:}
The case $|L_o(G)|=2$ has been proved by Theorems \ref{Wang} and \ref{Majie}. Now we suppose $|L_o(G)|\geq 3$.
We prove it by induction on $n=|V(G)|$. Obviously, it is true when $n\leq 2|L_o(G)|+1$. Assume that it is true for graphs with less than $n$ vertices.
If $G$ contains cut vertex, let $B$ be one of its blocks. Note that $|L_o(B)|\leq |L_o(G)|$. If $|L_o(B)|= |L_o(G)|$, 
then $w(B)\leq w(G)\leq 2|L_o(B)|$. Since $|V(B)|<n$, by the induction hypothesis, $\chi(B)\leq 2|L_o(B)|= 2|L_o(G)|$. If $|L_o(B)|< |L_o(G)|$, then
by Theorem \ref{AG}, $\chi(B)\leq 2|L_o(B)|+2\leq 2|L_o(G)|$. Hence, $\chi(G)\leq \max\{\chi(B)~|~B~is ~a~block~of~G\}\leq 2|L_o(G)|$.

Now we suppose that $G$ is 2-connected. If $\delta(G)\geq 2|L_o(G)|$, by  applying  Theorem \ref{2connected} on $k=|L_o(G)|$, we get $\chi(G)\leq 2|L_o(G)|$.  If $\delta(G)\leq 2|L_o(G)|-1$, let $v$ be the vertex with $d(v)\leq 2|L_o(G)|-1$ and let $G'=G-v$. Note that $|L_o(G')|\leq |L_o(G)|$. If $|L_o(G')|=|L_o(G)|$, then $w(G')\leq w(G)\leq 2|L_o(G')|$.
By the induction hypothesis, $\chi(G')\leq 2|L_o(G')|=2|L_o(G)|$. 
If $|L_o(G')|<|L_o(G)|$, by Theorem \ref{AG}, $\chi(G')\leq 2|L_o(G')|+2\leq 2|L_o(G)|$. In both cases, we get $\chi(G')\leq 2|L_o(G)|$. Since $d(v)\leq 2|L_o(G)|-1$, we see that $\chi(G)\leq 2|L_o(G)|$.\bqed

\vskip 2mm
\noindent{\bf{Proof of Theorem \ref{nok2L+2}}:}
We prove it by induction on $n=|V(G)|$. Obviously, it is true when $n\leq 2|L_e(G)|+2$. Assume that it is true for graphs with less than $n$ vertices.
If $G$ contains cut vertex, let $B$ be one of its blocks. Note that $|L_e(B)|\leq |L_e(G)|$. If $|L_e(B)|= |L_e(G)|$, 
then $w(B)\leq w(G)= 2|L_e(B)|+1$. Since $|V(B)|<n$, by the induction hypothesis, $\chi(B)\leq 2|L_e(B)|+1= 2|L_e(G)|+1$. If $|L_e(B)|< |L_e(G)|$, then
by Theorem \ref{MS}, $\chi(B)\leq 2|L_e(B)|+3\leq 2|L_e(G)|+1$. Hence, $\chi(G)\leq \max\{\chi(B)~|~B~is ~a~block~of~G\}\leq 2|L_e(G)|+1$.

Now we suppose that $G$ is 2-connected. If $\delta(G)\geq 2|L_e(G)|+1$, by  applying Theorem \ref{2connectedeven} on $k=|L_e(G)|$, we get $\chi(G)\leq 2|L_e(G)|+1$.  If $\delta(G)\leq 2|L_e(G)|$, let $v$ be the vertex with $d(v)\leq 2|L_e(G)|$ and let $G'=G-v$. Note that $|L_e(G')|\leq |L_e(G)|$. If $|L_e(G')|=|L_e(G)|$, then $w(G')\leq w(G)\leq 2|L_e(G')|+1$.
By the induction hypothesis, $\chi(G')\leq 2|L_e(G')|+1=2|L_e(G)|+1$. 
If $|L_e(G')|<|L_e(G)|$, by Theorem \ref{MS}, $\chi(G')\leq 2|L_e(G')|+3\leq 2|L_e(G)|+1$. In both cases, we get $\chi(G')\leq 2|L_e(G)|+1$. Since $d(v)\leq 2|L_e(G)|$, we see that $\chi(G)\leq 2|L_e(G)|+1$.\bqed

\section{Proof of Theorem \ref{2connected}}

The first part of Theorem \ref{2connected} has been proved by Theorem \ref{LWZ}. Now we suppose $|L_o(G)|=k$. By Theorem \ref{LWZ}, these $k$ odd cycles have consecutive odd lengths.
We remain to prove that $2|L_o(G)|+1=2k+1=\ell_o(G)$, and either $K_{2k+1}\ss G$ or $\chi(G)\leq 2k$.
Let $C:=v_0v_1\cdots v_{\ell-1}v_0$ be the cycle with $|C|=\ell_o(G)=\ell$. We say that the sequence $v_0,v_1,\ldots,v_{\ell-1},v_0$ is the \emph{forward} direction of $C$, and denote as $\overrightarrow{C}$. And the sequence $v_0,v_{\ell-1},\ldots,v_{1},v_0$ is the \emph{backward} direction of $C$, and denote as $\overleftarrow{C}$.
For convenience, every subscript $i$ of $v_i$ is taken module $\ell$.

\vskip 2mm

\noindent{\bf{Case 1.}}  $V(C)=V(G)$. For every $x\in V(G)$,
 since $d(x)\geq 2k$, there are at least $2k-2$ chords  incident with $x$. Denote $2k-2$ of them as $xx_j$ with $j=1,\ldots, 2k-2$. 
 Without loss of generality, suppose $x_1,\ldots,x_{2k-2}$ follow each other on $C$ in the forward direction.

Observe that each chord $xx_j$ separates the odd cycle $C$ into two small cycles $C_{1,j}=x\overrightarrow{C}x_jx$ and $C_{2,j}=x\overleftarrow{C}x_jx$. One is an even cycle and the other  is an odd cycle. 
It is easy to see that $C_{1,j}$ and $C_{2,j}$
are smaller than $C$ and  no two odd cycles $C_{1,j_1},C_{1,j_2}$ ($C_{2,j_1},C_{2,j_2}$) have the same length for any $j_1\not=j_2$. Since $|L_o(G)|=k$, 
 there are exactly $k-1$ odd cycles $C_{1,j}$ and $k-1$ 
 odd cycles $C_{2,j}$.
 Let $j_1<\cdots<j_{k-1}$ be the indices such that $C_{1,j_1},\ldots, C_{1,j_{k-1}}$ are odd. 
 And let $j_1'>\cdots>j_{k-1}'$ be the indices such that $C_{2,j_1'},\ldots, C_{2,j_{k-1}'}$ are odd.
 Then  $|x\overrightarrow{C}x_{j_1}|=|x\overleftarrow{C}x_{j_1'}|=e_1$ and $e_1$ is even.

Since $C_{1,j_1},\ldots, C_{1,j_{k-1}},C$  have consecutive odd lengths, these $k$ odd cycles have lengths $e_1+1,e_1+3,\ldots,e_1+2k-1$, implying that $|x_{j_{k-1}}\overrightarrow{C}x|=3$. If $x_{j_{k-1}}=x_{2k-2}$, then we have an odd cycle $xx_{2k-2}\overleftarrow{C}x_{j_1'}x$ smaller than $C_{2,j_1'}$, a contradiction to $|L_o(G)|=k$. Hence, $x_{j_1'}=x_{2k-2}$. And $x_{j_{k-1}},x_{j_1'},x$ follow each other on $C$ in the forward direction. 
Since $|x\overleftarrow{C}x_{j_1'}|=e_1$ is even and $|x_{j_{k-1}}\overrightarrow{C}x|=3$, we get $e_1=2$. Then $2|L_o(G)|+1=\ell_o(G)=2k+1$ and $|V(G)|=|V(C)|=2k+1$.
 Since $\delta(G)\geq 2k$, $G\cong K_{2k+1}$.

\vskip 3mm

Now we suppose $V(C)\not=V(G)$. We choose a longest odd cycle $C$ in $G$ such that: 

(1) if possible, $G-V(C)$ is not an independent set;

(2) if $G-V(C)$ is an independent set, $d(v)$ is as small as possible for every $v\in V(G)\backslash V(C)$.

\vskip 1mm

\noindent{\bf{Case 2.}}  $G-V(C)$ is an independent set. Let $x\in V(G)\backslash V(C)$.
Since $d(x)\geq 2k$, denote $2k$ neighbors of $x$ in $C$ as $x_0,\ldots,x_{2k-1}$, where they follow each other on $C$ in  the forward direction.
Note that for every $j\not=0$, the path $x_0xx_j$ separates $C$ into two cycles $C_{1,j}=x_0\overrightarrow{C}x_jxx_0$ and $C_{2,j}=x_0\overleftarrow{C}x_jxx_0$. One is an even cycle and the other is an odd cycle.  Since $|L_o(G)|= k$,  there are     $k$ odd $C_{1,j}$  or  $k$ odd $C_{2,j}$. 
Without loss of generality, suppose there are     $k$ odd $C_{1,j}$. Let $x_{j_1},\ldots,x_{j_k}$ be the neighbors of $x$ such that $x_0,x_{j_1},\ldots,x_{j_k}$ follow each other on $C$ in the forward direction and $C_{1,j_i}$ is an odd cycle  for every $1\leq i\leq k$. Then $C_{1,j_1},\ldots,C_{1,j_k}$ have consecutive odd lengths. And $|x_{j_i}\overrightarrow{C}x_{j_{i+1}}|=2$ for every $1\leq i\leq k-1$, $|C_{1,j_k}|=|C|$ and $|x_{j_{k}}\overrightarrow{C}x_0|=2$.
 Let $|x_0\overrightarrow{C}x_{j_1}|=O_1$. Then $O_1$ is odd and we have $k$ odd cycle lengths $O_1+2,O_1+4,\ldots,O_1+2k$.
Since $O_1+2$ is the smallest odd cycle length,
there is no neighbor of $x$  inside $x_0\overrightarrow{C}x_{j_1}$. Since $d(x)\geq 2k$, there must exist  $x^*\in N(x)$ that lies inside $x_{j_i}\overrightarrow{C}x_{j_{i+1}}$ for some $i$.
Then $xx^*x_{j_i}x$ is a triangle, implying that $O_1=1$ and $2|L_o(G)|+1=\ell_o(G)=2k+1$.

Now
we remain to prove that either $K_{2k+1}\ss G$ or $\chi(G)\leq 2k$. For convenience, suppose  $w(G)\leq 2k$, and we will prove $\chi(G)\leq 2k$. We first show that $d(x)=2k$ for any $x \in V(G)\backslash V(C)$. Note that $|C|=2k+1$, $G-V(C)$ is an independent set and  $d(x)\geq 2k$. Since $G[V(C)]\not\cong K_{2k+1}$, there exists an vertex $v_i\in V(C)$ and $d_C(v_i)\leq 2k-1$. If  $d(x)= 2k+1$, we can construct another longest odd cycle $C'=xv_{i+1}\overrightarrow{C}v_{i-1}x$. By our choice (1) of the longest odd cycle, $G-V(C')$ is also an independent set.  But then $v_i\in G-V(C')$ and  $d(v_i)=d_C(v_i)+1\leq 2k<d(x)$, a contradiction to our choice (2) of $C$.
Hence, $d(x)=2k$ for any $x \in V(G)\backslash V(C)$.

 Since $G[V(C)]\not\cong K_{2k+1}$, $\chi(G[V(C)])\leq 2k$. If $\chi(G[V(C)])\leq 2k-1$, then each vertex out of $C$ can choose the $2k^{th}$ color, implying that $\chi(G)\leq 2k$. 
 If $\chi(G[V(C)])=2k$,
note that $|V(C)|=2k+1$, then $G[V(C)]$ must contain $K_{2k}$. Denote the remaining vertex of $V(C)\backslash V(K_{2k})$ as $v$ and $v_{i_1}\in N(v)\cap V(K_{2k})$.
If $|V(G)\backslash V(C)|\geq 2$, denote two of them by $u_1,u_2$. Since $d(u_i)=2k$ and $w(G)\leq 2k$, $v\in N(u_i)$ and $|N(u_i)\cap K_{2k}|=2k-1$. Since $2k-1\geq 5$, we can find a neighbor $v_{i_2}$ of $u_1$ in  $V(K_{2k})-v_{i_1}$, and two neighbors $v_{i_3},v_{i_4}$ of $u_2$ in 
$V(K_{2k})\backslash\{v_{i_1},v_{i_2}\}$. But we can find a $(2k+3)$-cycle  in $K_{2k}\cup v_{i_1}vu_1v_{i_2}\cup v_{i_3}u_2v_{i_4}$, a contradiction.

Now we suppose $V(G)\backslash V(C)=\{x\}$. Since $w(G)= 2k$ and $d(x)=2k$, $v\in N(x)$ and there is only one vertex $x'$ in $ V(K_{2k})\backslash N(x)$. Then there is also only one vertex $v'$ in $ V(K_{2k})\backslash N(v)$. If $x'=v'$, then $V(G)-x'$ induces a clique of size $2k+1$, a contradiction. Hence, $G\cong K_{2k+2}^{2-}$,   By Lemma \ref{2-}, $\chi(G)\leq 2k$.

\vskip 3mm

Now we remain to consider that there exists a longest odd cycle $C$ such that $G-V(C)$ is not an independent set.  We will prove that this cannot happen by the assumption $|L_o(G)|= k$ and $C$ is the longest odd cycle.
 We choose a longest path $P$ in $G-V(C)$ with two endpoints $A$ and $B$ such that 

(1) if possible, $N(A)\cap V(C)\not=\emptyset$ and $N(B)\cap V(C)\not=\emptyset$;
 
 (2) the furthest neighbor of $A$ in $P$ is as close to $B$ in $P$ as possible;

 (3) the furthest neighbor of $B$ in $P$ is as close to $A$ in $P$ as possible.

\vskip 2mm

 \noindent{\bf{Case 3.}}
$N(A)\cap C=\emptyset$ or $N(B)\cap C=\emptyset$. 
 Without loss of generality, suppose $N(A)\cap C=\emptyset$.
Then all the neighbors of $A$ are in $P$. Denote the furthest neighbor of $A$ in $P$ as $A_0$, and the predecessor of $A_0$ in $P$ as $A'$. Then $A'PAA_0PB$ is also the longest path in $G-V(C)$. By the choices (1) and (2) of $P$, $N(A')\ss APA_0$ and $A_0$ is also the  furthest neighbor of $A'$ in $P$. Let $C'=APA'A_0A$ and $x\in \{A,A'\}$. Then     $N(x)\ss C'$,  and $AA_0A'\ss C'$.

Since $d(x)\geq 2k$, there are $2k-2$ chords $x_1x,\ldots,x_{2k-2}x$ in $C'$. 
Let $H_x$ be the union of  $C'$
and these $2k-2$ chords. If $H_x$ is non-bipartite, by applying Lemma \ref{non-bipartite} on $(H,t)=(H_x,k)$, then $|L(H_x)|\geq k$. Together with $C$, by Lemma \ref{C'}, we get $|L_o(G)|\geq k+1$, a contradiction.

If $H_x$ is bipartite for any $x\in \{A,A'\}$, then let $H=H_A\cup H_{A'}$. It is easy to check that $AA_0A'\ss C'$ and $H$ is bipartite.
Since $G$ is 2-connected, we can find two disjoint paths $s_1P_1t_1$ and $s_2P_2t_2$ connecting $V(C)$ and $V(C')$, where $s_1,s_2\in V(C)$ and $t_1,t_2\in V(C')$.
Applying Lemma \ref{bipartite} on $(H,y,z)=(H,t_1,t_2)$, we can find $k+1$ ($t_1,t_2$)-paths $Q_1,\ldots,Q_{k+1}$ of different lengths but of the same parity lengths in $G[V(C')]$. Since $C$ is an odd cycle, we find either  $k+1$ odd cycles $s_1\overrightarrow{C}s_2P_2t_2Q_it_1P_1s_1$ of different lengths, or $k+1$ odd cycles $s_1\overleftarrow{C}s_2P_2t_2Q_it_1P_1s_1$ of different lengths, a contradiction.


\vskip 2mm

\vskip 3mm

Now we suppose that  $N(A)\cap C\not=\emptyset$ and $N(B)\cap C\not=\emptyset$. 
Without loss of generality, suppose $|N(A)\cap V(P)|\geq |N(B)\cap V(P)|$.
 Let $|N(A)\cap V(C)|=p\geq 1$ and $|N(A)\cap V(P)|= q+1\geq 1$.  Since $d(A)\geq 2k$, $p+q+1\geq 2k$. We will use the  following claim to prove the remaining cases.

\begin{claim}\label{pq}(Gy\'arf\'as \cite{G})
    Assume that $A$ is adjacent to $x_1,\ldots,x_p$ in  $C$ and $q+1$ vertices in $P$. If $B$ is adjacent to $y\in V(C)\backslash \{x_1,\ldots,x_p\}$, then $|L_o(G)|\geq \lceil \frac{p}{2} \rceil+q$.
\end{claim}

\noindent{\bf{Case 4.}}
$N(A)\cap C\not=N(B)\cap C$. 
Let $y\in N(B)\cap V(C)\backslash N(A)$. By Claim \ref{pq}, $|L_o(G)|\geq \lceil \frac{p}{2} \rceil+q\geq \lceil\frac{2k-q-1}{2} \rceil+q=k+\lceil \frac{q-1}{2} \rceil$. Since $|L_o(G)|= k$, we get $q=0,1$. 
\vskip 2mm

\noindent{\bf{Subcase 4.1.}}
$q=0$. Since $|N(A)\cap V(C)|\geq 2k-1$, denote $x_1,\ldots,x_{2k-1}$
 as its $2k-1$ neighbors in $C$, and $y, x_1,\ldots,x_{2k-1}$ 
 follow each other on $C$ in the forward direction. Let $|P|=s>0$.

Note that the path $x_jAPBy$ separates $C$ into two cycles $C_{1,j}=APBy\overrightarrow{C}x_jA$ and $C_{2,j}=APBy\overleftarrow{C}x_jA$. One is an odd cycle and the other is an even cycle. Since $|L_o(G)|= k$, we have $k$ odd cycles $C_{1,j}$ or $k$ odd cycles $C_{2,j}$.
Without loss of generality, let $j_1<\cdots<j_{k}$ be the indices such that $C_{1,j_1},\ldots, C_{1,j_{k}}$ are odd. Then $|C|=|C_{1,j_{k}}|$ and  $|x_{j_k}\overrightarrow{C}y|=s+2$. Let $|y\overrightarrow{C}x_{j_1}|=a>0$. Since $C_{1,j_1},\ldots, C_{1,j_{k}}$ have consecutive odd lengths, they have lengths $s+a+2,\ldots,s+a+2k$, respectively.

 Suppose that there exists $y'\in N(B)\cap V(x_{j_k}\overrightarrow{C}x_{j_1})\backslash\{y,x_{j_k}\}$. If $y,y'$ follow each other on $C$ in the forward direction and $|y\overrightarrow{C}y'|$ is even or odd,  then $APBy'\overrightarrow{C}x_{j_1}A$ or $By\overrightarrow{C}y'B$ is an odd cycle smaller than $C_{1,j_1}$, a contradiction. 
If $y',y$ follow each other on $C$ in the forward direction and  $|y'\overrightarrow{C}y|$ is even, then $APBy'\overrightarrow{C}x_{j_k}A$ is an odd cycle larger than $C_{1,j_k}$, a contradiction.
Now we suppose that  $y',y$ follow each other on $C$ in the forward direction  and $|y\overrightarrow{C}y'|$ is odd.
Since $|x_{j_k}\overrightarrow{C}y|=s+2$, we have $|y'\overrightarrow{C}y|\leq s+1$.
 Then the odd cycle $By'\overrightarrow{C}yB$ has length at most $s+3\leq s+a+2=|C_{1,j_1}|$, implying that $a=1$, $|y'\overrightarrow{C}y|=s+1$ and  $|x_{j_k}\overrightarrow{C}y'|=1$.

The above arguments state that we have only two possible neighbors $y,y'$ of $B$ in $V(x_{j_k}\overrightarrow{C}x_{j_1})-x_{j_k}$. Since $1=|N(A)\cap V(P)|\geq |N(B)\cap V(P)|$,
we have  $|N(B)\cap V(C)|\geq 2k-1$. Then there are at least $2k-3$ neighbors of $B$ in $V(x_{j_1}\overrightarrow{C}x_{j_k})-x_{j_1}$. Since  
$|V(x_{j_1}\overrightarrow{C}x_{j_k})-x_{j_1}|=2k-2$,  there must exist a 3-cycle in a subgraph induced by $V(C)\cup \{B\}$. But $|C_{1,j_1}|=s+a+2>3$, a contradiction.

\vskip 2mm

\noindent{\bf{Subcase 4.2.}}
$q=1$. Since $|N(A)\cap V(C)|\geq 2k-2$, denote $x_1,\ldots,x_{2k-2}$
 as its $2k-2$ neighbors in $C$, and $y, x_1,\ldots,x_{2k-2}$ 
 follow each other on $C$ in the forward direction.
We first prove the following claim.

\begin{claim}\label{q=1}
   There exist no two 
$(A,B)$-paths $P_1,P_2$ in $G[V(P)]$ such that $|P_1|\not=|P_2|$ and $|P_1|\equiv |P_2| \pmod{2}$. 
\end{claim}

\noindent{\bf{Proof.}}
Suppose to the contrary that there are  two 
$(A,B)$-paths $P_1,P_2$ in $G[V(P)]$ such that
  $|P_1|=s+t$ and $|P_2|=t$, where $s>0$ is even and $t>0$.
Note that the path $x_jAP_iBy$ separates $C$ into two cycles $C_{1,j}^i=AP_iBy\overrightarrow{C}x_jA$ and $C_{2,j}^i=AP_iBy\overleftarrow{C}x_jA$. One is an odd cycle and the other is an even cycle. If there exists $k$ odd cycles $C_{1,j}^1$ or $k$ odd cycles $C_{2,j}^1$, let $C_{1,j^*}^1$ or $C_{2,j^*}^1$ be the smallest one, respectively.
    Together with $C_{1,j^*}^2$ or $C_{2,j^*}^2$, we get $k+1$ odd cycles of different lengths, a contradiction. Hence, we have 
$k-1$ odd cycles $C_{1,j}^i$ and $k-1$ odd cycles $C_{2,j}^i$ for $i\in \{1,2\}$. Let $j_1<\cdots<j_{k-1}$ be the indices such that $C_{1,j_1}^i,\ldots, C_{1,j_{k-1}}^i$ are odd.

 Let $|y\overrightarrow{C}x_{j_1}|=a$. Then $|C_{1,j_1}^2|=t+a+2$.
Since there are  $k$ odd cycles $C_{1,j_1}^2,C_{1,j_1}^1,\ldots, $ $C_{1,j_{k-1}}^1$ of different lengths and they have consecutive odd lengths, we can check that these cycles have lengths 
$t+a+2,t+a+4,\ldots, t+a+2k$,  $s=2$ and $|x_{j_i}\overrightarrow{C}x_{j_{i+1}}|=2$.  Since $|C_{1,j_{k-1}}^1|=|C|$,  $|x_{j_{k-1}}\overrightarrow{C}y|=t+4$.

Let $y'\in N(B)\cap V(C)-y$.
Since $2=|N(A)\cap V(P)|\geq |N(B)\cap V(P)|$, $|N(B)\cap V(C)|\geq 2k-2$. If $y'=x_{j_1}$, then one of $\{AP_2Bx_{j_1}A,By\overrightarrow{C}x_{j_1}B\}$ is an odd cycle smaller than $C_{1,j_1}^2$, a contradiction. 
Since  $|V(x_{j_1}\overrightarrow{C}x_{j_{k-1}})\backslash \{x_{j_1},x_{j_{k-1}}\}|=2k-5$, we have at least $3$ neighbors of $B$  in $V(x_{j_{k-1}}\overrightarrow{C}x_{j_1})-x_{j_1}$.

 Suppose  $y'\in  V(x_{j_{k-1}}\overrightarrow{C}x_{j_1})\backslash \{x_{j_1},x_{j_{k-1}}\}$. If $y,y'$ follow each other on $C$ in the forward direction and $|y\overrightarrow{C}y'|$ is even or odd,  then $AP_2By'\overrightarrow{C}x_{j_1}A$ or $By\overrightarrow{C}y'B$ is an odd cycle smaller than $C_{1,j_1}^2$, a contradiction. 
If $y',y$ follow each other on $C$ in the forward direction  and $|y'\overrightarrow{C}y|$ is even, then $AP_1By'\overrightarrow{C}x_{j_{k-1}}A$ is an odd cycle larger than $C_{1,j_{k-1}}^1$, a contradiction.

Now we suppose that  $y',y$ follow each other on $C$ in the forward direction  and $|y'\overrightarrow{C}y|$ is odd.
Let $|y'\overrightarrow{C}y|=a'$.
Since $Ax_{j_{k-1}}\overrightarrow{C}y'BP_2A$ is an odd cycle of length $2t+6-a'$ and  $By'\overrightarrow{C}yB$ is an odd cycle of length $a'+2$,  one of them has length at most $ t+4<t+a+4=|C_{1,j_1}^1|$. That means one of them has length $t+a+2=|C_{1,j_1}^2|\leq t+4$.
And  we get $a=1$ or $2$.
 If $a=1$, then  $\{|Ax_{j_{k-1}}\overrightarrow{C}y'BP_2A|,|By'\overrightarrow{C}yB|\}=\{t+3,t+5\}$. That means $a'=t+1$ or $t+3$. If $a'=t+1$, then $|y'\overrightarrow{C}x_{j_1}|=t+2$. Since $|x_{j_1}AP_1By'|=t+4$, we get an odd cycle $x_{j_1}AP_1By'\overleftarrow{C}x_{j_1}$ larger than $C$, a contradiction. Hence, $a'=t+3$. 
If $a=2$, then both $Ax_{j_{k-1}}\overrightarrow{C}y'BP_2A$
and  $By'\overrightarrow{C}yB$ have length $t+4$, $a'=t+2$ and  $|x_{j_{k-1}}\overrightarrow{C}y'|=2$. 

In  both cases of $a$,  $B$ has at most three neighbors $y,y',x_{j_{k-1}}$ in $V(x_{j_{k-1}}\overrightarrow{C}x_{j_1}) -x_{j_1}$. 
Hence, $N(B)\cap V(C)=\{y,y'\}\cup V(x_{j_1}\overrightarrow{C}x_{j_{k-1}})-x_{j_1}$. But there is an odd cycle $Bx_{j_{2}}x_{j_{2}}^-B$ of length $3<t+a+2$, a contradiction. 
\bqed

\vskip 2mm

Recall that the furthest neighbor of $A$ is $A_0$  and $A_0$ may be $B$. If $|APA_0|$ is odd, then  $APB$ and $AA_0PB$ have different but the same parity lengths, a contradiction to Claim \ref{q=1}. 
 Let $|APA_0|=s_1$ and $|A_0PB|=s_2$. Then $s_1\geq 2$ is even, $APA_0A$ is an odd cycle of length $s_1+1$, and $s_2\geq 0$. 
Note that the path $x_jAPBy$ separates $C$ into two cycles $C_{1,j}=APBy\overrightarrow{C}x_jA$ and $C_{2,j}=APBy\overleftarrow{C}x_jA$. One is an odd cycle and the other is an even cycle. 
Since $|L_o(G)|= k$ and $APA_0A$ is an odd cycle, we have 
$k-1$ odd cycles $C_{1,j}$ and $k-1$ odd cycles $C_{2,j}$. Let $j_1<\cdots<j_{k-1}$ be the indices such that $C_{1,j_1},\ldots, C_{1,j_{k-1}}$ are odd.
Then 
we have $k$ odd cycles $APA_0A,C_{1,j_1},\ldots,C_{1,j_{k-1}}$. 
Since $APA_0A,C_{1,j_1},\ldots,C_{1,j_{k-1}}$ have consecutive odd lengths,  they have lengths $s_1+1,s_1+3,\ldots,s_1+2k-1$, respectively. 
Then $s_2=0$, $AB\in E(G)$ and  $|y\overrightarrow{C}x_{j_1}|=1$. But $Ax_{j_1}\overrightarrow{C}yBA$ is an odd cycle larger than $C$, a contradiction.

\vskip 3mm

\noindent{\bf{Case 5.}}
$N(A)\cap C=N(B)\cap C$. Recall that $|N(A)\cap V(C)|=p\geq 1$, $|N(A)\cap V(P)|= q+1\geq 1$ and $p+q+1\geq 2k$. 
If $p=1$, denote $\{x_0\}=N(A)\cap C=N(B)\cap C$.  
Let $C'=APBx_0A$ and $x\in\{A,B\}$. Then  $N(x)\ss C'$ and there are $2k-2$ chords $x_1x,\ldots,x_{2k-2}x$ in $C'$. Let $H_x$ be the union of  $C'$
and $2k-2$ chords $xx_j$. If $H_x$ is non-bipartite, by applying Lemma \ref{non-bipartite} on $(H,t)=(H_x,k)$, then $|L(H_x)|\geq k$. Together with $C$, by Lemma \ref{C'}, we get $|L_o(G)|\geq k+1$, a contradiction. If $H_x$ is bipartite for any $x\in\{A,B\}$, let $H$ be the union of $C'$ and $2k-2$ chords incident with  $A$ and  $B$. Since $G$ is 2-connected, we can find a path $s_1P_1t_1$  connecting $V(C)$ and $V(C')$, where $s_1\in V(C)$ and $t_1\in V(C')$.
Applying  Lemma \ref{bipartite} on $(H,y,z)=(H,x_0,t_1)$, we can find $k+1$ ($x_0,t_1$)-paths $Q_1,\ldots,Q_{k+1}$ of different lengths but of the same parity lengths in $G[V(C')]$.
Since $C$ is an odd cycle, either we find $k+1$ odd cycles $s_1\overrightarrow{C}x_0Q_it_1P_1s_1$ of different lengths, or $k+1$ odd cycles $s_1\overleftarrow{C}x_0Q_it_1P_1s_1$ of different lengths, a contradiction.

If $p\geq 2$, by applying Claim \ref{pq} with $p-1$ in the role of $p$, we have $|L_o(G)|\geq \lceil \frac{p-1}{2}\rceil +q\geq \lceil \frac{2k-q-2}{2}\rceil +q=k-1+ \lceil \frac{q}{2}\rceil$. Since $|L_o(G)|= k$, $q\leq 2$. Let $|P|=s>0$.

\vskip 2mm

\noindent{\bf{Subcase 5.1.}} $|N(A)\cap C|=|N(B)\cap C|\geq 2k-2$. Denote  
 $x,x_1,\ldots,x_{2k-3}$
 as their $2k-2$ neighbors in $C$.
 For convenience, let $ x,x_1,\ldots,x_{2k-3}$ 
 follow each other on $C$ in the forward direction.
Note that the path $x_jAPBx$ separate $C$ into two cycles $C_{1,j}=APBx\overrightarrow{C}x_jA$ and $C_{2,j}=APBx\overleftarrow{C}x_jA$. One is an odd cycle and the other is an even cycle. 
It is easy to see that at least $k-1$ $C_{1,j}$ are odd or at least $k-1$ $C_{2,j}$ are odd.
Without loss of generality, let $j_1<\cdots<j_{k-1}$ be the indices such that $C_{1,j_1},\ldots, C_{1,j_{k-1}}$ are odd. 
Let $|P|=s$. If $s$ is even,  then we have an odd cycle $Ax\overrightarrow{C}x_{j_1}A$ smaller than $C_{1,j_1}$, implying that we have $k$ odd cycles $Ax\overrightarrow{C}x_{j_1}A,C_{1,j_1},\ldots, C_{1,j_{k-1}}$ of consecutive odd lengths. Then $|x_{j_1}\overrightarrow{C}x_{j_{2}}|=2$. But we get an odd cycle $x_{j_1}APBx_{j_2}\overrightarrow{C}x_{j_1}$ larger than $C$, a contradiction.

If $s$ is odd, then 
 there are $k$ odd cycles $APBxA,C_{1,j_1},\ldots,C_{1,j_{k-1}}$ of consecutive odd lengths $s+2,s+4,\ldots,s+2k$, respectively. 
If $x_{j_i}^-$ is a common neighbor of $A$ and $B$, then we get a path $x_{j_i}APBx_{j_i}^-$ of length $s+2\equiv 1$. And we have an odd cycle $x_{j_i}APBx_{j_i}^-\overleftarrow{C}x_{j_i}$ larger than $C$, a contradiction.
Suppose $A$ and $B$ have a common neighbor $x'\in V(x_{j_{k-1}}\overrightarrow{C}x)\backslash \{x,x_{j_{k-1}}\}$. Since $|x_{j_{k-1}}\overrightarrow{C}x|=s+2$ is odd,  either $x_{j_{k-1}}\overrightarrow{C}x'$ or $x'\overrightarrow{C}x$ is a path of odd length less than $s+2$. Since both two paths $x_{j_{k-1}}APBx'$ and $x'APBx$ have odd length $s+2$,  we have an odd cycle $x_{j_{k-1}}APBx'\overrightarrow{C}x_{j_{k-1}}$ or $x'APBx\overrightarrow{C}x'$ larger than $C$, a contradiction.
The above arguments state that $N(A)\cap C=N(B)\cap C=\{x,x_{j_1},\ldots,x_{j_{k-1}}\}$, a contradiction to $|N(A)\cap C|=|N(B)\cap C|\geq 2k-2$.

\vskip 2mm

\noindent{\bf{Subcase 5.2.}} $|N(A)\cap C|=|N(B)\cap C|\leq 2k-3$. Since $q+1\leq 3$, we get 
$|N(A)\cap C|=|N(B)\cap C|= 2k-3$ and $|N(A)\cap P|=|N(B)\cap P|= 3$. Denote $2k-3$ vertices of $N(A)\cap C$ as $x,x_1,\ldots,x_{2k-4}$. And without loss of generality, $x,x_1,\ldots,x_{2k-4}$ follow each other on $C$ in the forward direction.
We first prove the following claim.

\begin{claim}\label{no_odd}
There exist no two $(A,B)$-paths $P_1,P_2\in G[V(P)]$ such that $|P_1|\equiv |P_2|\equiv 1 $ and $|P_1|\not=|P_2|$.     
\end{claim}

\noindent{\bf{Proof.}}
Suppose to the contrary that there exist two $(A,B)$-paths $P_1,P_2\ss G[V(P)]$ such that
$|P_1|=s_1$, $|P_2|=s_2$, and $ s_1>s_2\geq 1$ are odd.
Note that the path $x_jAP_1Bx$ separates $C$ into two cycles $C_{1,j}=AP_1Bx\overrightarrow{C}x_jA$ and $C_{2,j}=AP_1Bx\overleftarrow{C}x_jA$. One is an odd cycle and the other is an even cycle. Since $AP_2BxA$ and $AP_1BxA$ are two odd cycles smaller than $C_{1,j}$ and $C_{2,j}$, we can check that 
 we have exactly
$k-2$ odd cycles $C_{1,j}$ and $k-2$ odd cycles $C_{2,j}$. Let $j_1<\cdots<j_{k-2}$ be the indices such that $C_{1,j_1},\ldots, C_{1,j_{k-2}}$ are odd. And let  $j_1'>\cdots>j_{k-2}'$ be the indices such that $C_{2,j_1'},\ldots, C_{2,j_{k-2}'}$ are odd.
Then we have $k$ odd cycles
$AP_2BxA,AP_1BxA, C_{1,j_1},\ldots, C_{1,j_{k-2}}$ and they have consecutive odd lengths, implying that they have lengths $s_2+2,\ldots,s_2+2k$,  $s_1=s_2+2$ and $|x_{j_{k-2}}\overrightarrow{C}x|=s_2+4$. We also have $k$ odd cycles
$AP_2BxA, AP_1BxA, C_{2,j_1'},\ldots, C_{2,j_{k-2}'}$ of lengths $s_2+2,\ldots,s_2+2k$. Then $|x\overleftarrow{C}x_{j_1'}|=2$ and $|x_{j_{k-2}}\overrightarrow{C}x_{j_1'}|=s_2+2$.
But we obtain an odd cycle $x_{j_{k-2}}AP_1Bx_{j_1'}\overrightarrow{C}x_{j_{k-2}}$ larger than $C$, a contradiction. \bqed

\vskip 2mm
Now we continue to prove.
Let $y$ be a neighbor of $A$ or $B$ in $P$, and $Ay\notin P$ or $By\notin P$, respectively. 
We say    $y$ is \emph{odd} if  $|APy|$ or $|BPy|$ is odd. Otherwise, $y$ is \emph{even}.

If $|P|$ is odd, by Claim \ref{no_odd}, $A$ and $B$ both have no odd neighbor. Let $y_A,y_B$ be an even neighbor of $A,B$, respectively. 
If $P=APy_APy_BPB$, then we get two  $(A,B)$-paths $APB$ and $Ay_APy_BB$ of different odd lengths, a contradiction to Claim \ref{no_odd}.
If $P=APy_BPy_APB$, denote $|APy_B|=s_1$, $|y_BPy_A|=s_2$ and $|y_APB|=s_3$. Since $|P|$ is odd and $y_A$ and $y_B$ are even, we get $s_i$ is odd for every $i\in \{1,2,3\}$.
Since $APB$ and $Ay_APy_BB$ are two paths of odd lengths $s_1+s_2+s_3$ and $s_2+2$, by Claim \ref{no_odd}, we get $s_1=s_3=1$. Note that $|N(A)\cap P|= 3$. Denote $y_A'$ as another even neighbor of $A$ in $P$. Since $s_1=s_3=1$, $P=Ay_BPy_A'Py_AB$. Then
$|y_BPy_A'|=s_2'<s_2$ and is odd. And we get two  paths $APB$ and $Ay_A'Py_BB$ of  odd lengths $s_2+2$ and $s_2'+2$, a contradiction to Claim \ref{no_odd}.

Now we suppose that $|P|$ is even. If $A$  has two even neighbors $y_1,y_2$,   then $Ay_1PB$ and $Ay_2PB$ are two paths of different odd lengths, a contradiction to Claim \ref{no_odd}. Then $A$ has at least one odd neighbor $y_A$.  Let $P_1=P$, $P_2=Ay_APB$, $|P_1|=s_1>0$ and $|P_2|=s_2>0$. Then $s_1>s_2$ and $s_1\equiv s_2\equiv 0 $.
Note that the path $x_jAP_iBx$ separates $C$ into two cycles $C_{1,j}^i=AP_iBx\overrightarrow{C}x_jA$ and $C_{2,j}^i=AP_iBx\overleftarrow{C}x_jA$ for every $i\in \{1,2\}$. One is an odd cycle and the other is an even cycle. If there exists $k-1$ odd cycles $C_{1,j}^2$ or $k-1$ odd cycles $C_{2,j}^2$, let $C_{1,j^{*}}^2$ or $C_{2,j^{*}}^2$ be the smallest one and  let $C_{1,j^{**}}^2$ or $C_{2,j^{**}}^2$ be the largest one, respectively.
 Together with $Ax\overrightarrow{C}x_{j^{*}}A$ and  $C_{1,j^{**}}^1$ or $Ax\overleftarrow{C}x_{j^{*}}A$ and $C_{2,j^{**}}^1$, we get $k+1$ odd cycles of different lengths, a contradiction. Hence, we have exactly
$k-2$ odd cycles $C_{1,j}^i$ and $k-2$ odd cycles $C_{2,j}^i$ for $i\in \{1,2\}$. Let $j_1<\cdots<j_{k-2}$ be the indices such that $C_{1,j_1}^i,\ldots, C_{1,j_{k-2}}^i$ are odd. And let $j_1'>\cdots>j_{k-2}'$ be the indices such that $C_{1,{j_1}'}^i,\ldots, C_{1,{j_{k-2}'}}^i$ are odd.
  Let $|x\overrightarrow{C}x_{j_1}|=a$. Then $a$ is odd.

Since there are $k$ odd cycles 
$Ax\overrightarrow{C}x_{j_1}A$, $C_{1,j_1}^2,\ldots, C_{1,j_{k-2}}^2,  C_{1,j_{k-2}}^1$ and they have consecutive odd lengths, they have lengths $a+2,\ldots,a+2k$. Then we can check that $s_2=2$ and $s_1=4$.
If $k\geq 4$, then  $|x_{j_1}\overrightarrow{C}x_{j_2}|=2$. And $x_{j_1}AP_1Bx_{j_2}\overrightarrow{C}x_{j_1}$ is  an odd cycle larger than $C$, a contradiction.
If $k=3$, then there are $3$ odd cycles 
$Ax\overrightarrow{C}x_{j_1}A$, $C_{1,j_1}^2,  C_{1,j_{1}}^1$ of lengths $a+2,a+4,a+6$.
And $|x_{j_1}\overrightarrow{C}x|=6$.
Similarly, there are $3$ odd cycles 
$Ax\overleftarrow{C}x_{j_1'}A$, $C_{1,j_1'}^2,  C_{1,j_{1}'}^1$ of lengths $a+2,a+4,a+6$. 
Then $|x\overleftarrow{C}x_{j_1'}|=a$ and $|x_{j_1'}\overleftarrow{C}x|=6$.
If $x_{j_1}=x_2$ and $x_{j_1'}=x_1$, then $a>6$ and $|x_1\overrightarrow{C}x_2|=a-6$. But $Ax_1\overrightarrow{C}x_2A$ is an odd cycle of length $a-4<a+2$, a contradiction.
Then $x_{j_1}=x_1$ and $x_{j_1'}=x_2$. 
Since $|x_2\overrightarrow{C}x|=|x\overrightarrow{C}x_1|=a$ is odd, $|x_1\overrightarrow{C}x_2|$ is odd. Since there are 3 odd cycles $Ax_1\overrightarrow{C}x_2A$, $AP_2Bx_1\overrightarrow{C}x_2A$ and $AP_1Bx_1\overrightarrow{C}x_2A$ of different lengths, we get $|x_1\overrightarrow{C}x_2|=a$. Then $|C|=3a=a+6$ and the three odd cycle lengths are $5,7,9$.
Recall that $P=P_1$, $P_2=Ay_APB$, $|P_1|=4$ and $|P_2|=2$. 
Denote the third neighbor of  $A$ in $P$ as $y_A'$. But any possible  position of $y_A'$ will  produce a triangle,
a contradiction to $|L_o(G)|=3$.
\bqed

\section{Proof of Theorem \ref{2connectedeven}}

The first part of Theorem \ref{2connectedeven} has been proved by Theorem \ref{LM}. Now we suppose $|L_e(G)|=k$. By Theorem \ref{LM}, these $k$ even cycles have consecutive even lengths.
We remain to prove that $2|L_e(G)|+2=2k+2=\ell_e(G)$, and either $K_{2k+2}\ss G$ or $\chi(G)\leq 2k+1$.
Let $C:=v_0v_1\cdots v_{\ell-1}v_0$ be a cycle with $|C|=\ell_e(G)=\ell$. We say that the sequence $v_0,v_1,\ldots,v_{\ell-1},v_0$ is the \emph{forward} direction of $C$, and denote as $\overrightarrow{C}$. And the sequence $v_0,v_{\ell-1},\ldots,v_{1},v_0$ is the \emph{backward} direction of $C$, and denote as $\overleftarrow{C}$.
For convenience, every subscript $i$ of $v_i$ is taken module $\ell$.

\vskip 2mm

\noindent{\bf{Case 1.}}  $V(C)=V(G)$. 
For every $x\in V(C)$,
 since $d(x)\geq 2k+1$, there are at least $2k-1$ chords  incident with $x$. 
Denote $2k-1$ of them  as $xx_1,\ldots,xx_{2k-1}$, and $x_1,\ldots,x_{2k-1}$ follow each other on $C$ in the forward direction. Let $H$ be the union of $C$ and $xx_1,\ldots,xx_{2k-1}$.
By  Lemma \ref{C2k-1},   $|V(G)|=|V(C)|=2k+2$. Since $\delta(G)\geq 2k+1$, $G\cong K_{2k+2}$.

\vskip 3mm

Now we suppose $V(C)\not=V(G)$. We choose a longest even cycle $C$ of $G$ such that: 

(1) if possible, $G-V(C)$ is not an independent set;

(2) if $G-V(C)$ is an independent set, $d(v)$ is as small as possible for every $v\in V(G)\backslash V(C)$.

\vskip 1mm

\noindent{\bf{Case 2.}}  $G-V(C)$ is an independent set. Let $x\in V(G)\backslash V(C)$.
Since $d(x)\geq 2k+1$, denote $2k+1$ neighbors of $x$ in $C$ as $x_0,\ldots,x_{2k}$, where they follow each other on $C$ in  the forward direction. For convenience, every subscript $i$ of $x_i$ is taken module $2k+1$.
For any vertex $x_i$ and $j\not=i$, the path $x_ixx_j$ separates $C$ into two cycles $C_{i,j}=x_i\overrightarrow{C}x_jxx_i$ and $C_{i,j}'=x_i\overleftarrow{C}x_jxx_i$. Both  are  even or odd cycles. It is easy to see that if $C_{j_1,j_2}$ and $C_{j_2,j_3}$ are even cycles, then $C_{j_1,j_3}$ is even for every distinct $j_1,j_2,j_3$. We say a subset $S\ss \{x_0,\ldots, x_{2k}\}$ is even if $C_{j_1,j_2}$ is even for any  $x_{j_1},x_{j_2}\in S$. Obviously, there must exist an even set with cardinality at least 2.

Let $S$ be the maximum even set.
We  claim that  $|S|\geq k+1$. Otherwise, $|S|\leq k$. Without loss of generality, suppose $x_0\in S$. Then $C_{0,j}$ is odd for any $x_j\notin S$. But  $\{x_0,\ldots, x_{2k}\}\backslash S$ is even with cardinality at least $k+1$, a contradiction.
Denote $k+1$ vertices of $S$ as $y_0,\ldots,y_{k}$, and $y_0,\ldots,y_{k}$ follow each other on $C$ in the forward direction. For every $y_i$, we have an even cycle $xy_i\overrightarrow{C}y_{i+j}x$ for $1\leq j\leq k$, where $i+j$ is taken module $k+1$. Then we have $k$ even cycles of different lengths. And $|xy_i\overrightarrow{C}y_{i+k}x|=|C|$, implying that $|y_{i+k}\overrightarrow{C}y_i|=2$. That is, $|y_{i-1}\overrightarrow{C}y_i|=2$. Then $|C|=2k+2$ and $2|L_e(G)|+2=\ell_e(G)=2k+2$.

Now
we remain to prove that either $K_{2k+2}\ss G$ or $\chi(G)\leq 2k+1$. For convenience, suppose  $w(G)\leq 2k+1$, and we will prove $\chi(G)\leq 2k+1$. We first prove that $d(x)=2k+1$ for any $x \in V(G)\backslash V(C)$. Note that $|C|=2k+2$, $G-V(C)$ is an independent set and  $d(x)\geq 2k+1$. Since $G[V(C)]\not\cong K_{2k+2}$, there exists a vertex $v_i\in V(C)$ and $d_C(v_i)\leq 2k$. If  $d(x)= 2k+2$, we can construct another longest even cycle $C'=xv_{i+1}\overrightarrow{C}v_{i-1}x$. By our choice (1) of the longest even cycle, $G-V(C')$ is also an independent set.  But then $v_i\in G-V(C')$ and  $d(v_i)=d_C(v_i)+1\leq 2k+1<d(x)$, a contradiction to our choice (2) of $C$.
Hence, $d(x)=2k+1$ for any $x \in V(G)\backslash V(C)$.

 Since $G[V(C)]\not\cong K_{2k+2}$, $\chi(G[V(C)])\leq 2k+1$. If $\chi(G[V(C)])\leq 2k$, then each vertex out of $C$ can choose the $(2k+1)^{th}$ color, implying that $\chi(G)\leq 2k+1$. 
 If $\chi(G[V(C)])=2k+1$,
note that $|V(C)|=2k+2$, then $G[V(C)]$ must contain $K_{2k+1}$. Denote the remaining vertex of $V(C)\backslash V(K_{2k+1})$ as $v$. Since $v$ has at least two neighbors in $V(K_{2k+1})$, denote two of them by $v_{i_1},v_{i_2}$.
If $|V(G)\backslash V(C)|\geq 2$, denote two of them by $u_1,u_2$. Since $d(u_i)=2k+1$ and $w(G)\leq 2k+1$, $v\in N(u_i)$ and  $|N(u_i)\cap K_{2k+1}|=2k$. Since $2k\geq 6$, we can find  two neighbors $v_{i_3},v_{i_4}$ and $v_{i_5},v_{i_6}$ of $u_1$ and $u_2$ in $V(K_{2k+1})\backslash\{v_{i_1},v_{i_2}\}$,  respectively. And $v_{i_1},v_{i_2},v_{i_3},v_{i_4},v_{i_5},v_{i_6}$ are distinct. Then we can find an even cycle of length $2k+4$ in $K_{2k+1}\cup v_{i_1}vv_{i_2}\cup v_{i_3}u_1v_{i_4}\cup v_{i_5}u_2v_{i_6}$, a contradiction.

Now we suppose $V(G)\backslash V(C)=\{x\}$. Since $w(G)= 2k+1$ and $d(x)=2k+1$, $v\in N(x)$ and there is only one vertex $x'$ in $ V(K_{2k+1})\backslash N(x)$. Then there is also only one vertex $v'$ in $ V(K_{2k+1})\backslash N(v)$. If $x'=v'$, then $V(G)-x'$ induces a clique of size $2k+2$, a contradiction. Hence, $G\cong K_{2k+3}^{2-}$.
By Lemma \ref{2-}, $\chi(G)\leq 2k+1$.

\vskip 3mm

Now we remain to consider that there exists a longest even cycle $C$ such that $G-V(C)$ is not an independent set. We will prove that this cannot happen by   our assumption $|L_e(G)|= k$ and $C$ is the largest even cycle. 
 We choose a longest path $P$ in $G-V(C)$ with two endpoints $A$ and $B$ such that 

(1) if possible, $N(A)\cap V(C)\not=\emptyset$ and $N(B)\cap V(C)\not=\emptyset$.
 
 (2) the furthest neighbor of $A$ in $P$ is as close to $B$ in $P$ as possible,

 (3) the furthest neighbor of $B$ in $P$ is as close to $A$ in $P$ as possible.

\vskip 2mm

 \noindent{\bf{Case 3.}}
$N(A)\cap C=\emptyset$ or $N(B)\cap C=\emptyset$. 
 Without loss of generality, let $N(A)\cap C=\emptyset$.
Then all the neighbors of $A$ are in $P$. Denote the furthest neighbor of $A$ in $P$ as $A_0$, and the predecessor of $A_0$ in $APB$ as $A'$. Then $A'PAA_0PB$ is also the longest path in $G-V(C)$. By the choices (1) and (2) of $P$, $N(A')\ss APA_0$ and $A_0$ is also the  furthest neighbor of $A'$ in $P$. Let $C^*=APA'A_0A$. Then we have the following claim.

\begin{claim}\label{A'e}
For any $x\in \{A,A'\}$,     $N(x)\ss V(C^*)$ and $AA_0A'\ss C^*$.
\end{claim}

We define the direction $APA'A_0A$ of $C^*$ as the \emph{forward} direction, and denote as $\overrightarrow{C^*}$. And define another direction of $C^*$ as the \emph{backward} direction, and denote as $\overleftarrow{C^*}$.
By Claim \ref{A'e}, there are $2k-1$ chords $x_jA$ in $C^*$ with $1\leq j\leq 2k-1$. 
 Without loss of generality, suppose $x_1,\ldots,x_{2k-1}$ follow each other on $C^*$ in the forward direction.

 \vskip 1mm
\noindent{\bf{Subcase 3.1.}} $C^*$ is even.  
Let $H$ be the union of $C^*$ and $x_1A,\ldots,x_{2k-1}A$.
By Lemma \ref{C2k-1},  $|C^*|=2k+2$.  
Since $d(A)\geq 2k+1$ and $N(A)\ss C^*$, $A$ is adjacent to every vertex in $V(C^*)-A$.
Let $x$ be any vertex in $V(C^*)\backslash\{A,A_0\}$. Then we can   
find another longest path $x\overleftarrow{C^*}Ax^+\overrightarrow{C^*}A_0PB$. By our choices (1) and (2) of $P$, $N(x)\ss C^*$. Then $x$ is adjacent to every vertex in $V(C^*)-x$ for every $x\in V(C^*)\backslash A_0$, implying that  $G[V(C^*)]\cong K_{2k+2}$.
Since $G$ is 2-connected, 
we can find two disjoint paths $s_1P_1t_1$ and $s_2P_2t_2$ connecting $V(C)$ and $V(C^*)$, where $s_1,s_2\in V(C)$ and $t_1,t_2\in V(C^*)$.
Let $Q_1$ be a $(s_1,s_2)$-path in $C$ of length at least $k+1$.
Since $G[V(C^*)]\cong K_{2k+2}$, we can find a   $(t_1,t_2)$-path $Q_2$ in $G[V(C^*)]$ of length at least $k+1$ and $|Q_2|\equiv |Q_1|+|P_1|+|P_2|$. Then we have an even cycle $s_1P_1t_1Q_2t_2P_2s_2Q_1s_1$ of length larger than $2k+2$, a contradiction.

\vskip 2mm

\noindent{\bf{Subcase 3.2.}} $C^*$ is odd. 
Let  $y_{j}A'$ be a chord of  $A'$ for $1\leq j\leq 2k-1$. Without loss of generality, suppose $y_{1},\ldots,y_{{2k-1}}$ follow each other on $C^*$ in the  forward direction.
Let $H_A$ be the union of $C^*$ and $x_1A,\ldots,x_{2k-1}A$, and $H_{A'}$ be the union of $C^*$ and $y_1A',\ldots,y_{2k-1}A'$.
Applying Lemma \ref{C2k-1odd} on $H\in \{H_A,H_{A'}\}$, we have $k$ even cycle lengths $4,6,\ldots,2k+2$ and  two possible situations for each $A$ and $A'$.
We will prove  $|C^*|=2k+3$ by two  cases.

(I) 
Suppose we have  $k$ even cycles $C_{A,j_1}',\ldots,C_{A,j_k}'$ and $k$ even cycles $C_{A',j_1},\ldots,C_{A',j_k}$.
By Lemma \ref{C2k-1odd}, we have  
$|A\overleftarrow{C^*}x_{j_{1}}|=3$ and $|x_{j_1}\overleftarrow{C^*}x_{j_2}|=2$.  Let $e_0=|y_{j_k}\overrightarrow{C^*}A'|$. Then $e_0$ is even.  If $e_0\geq 4$, then we can check that $y_{j_k},x_{j_2},x_{j_1},A',A$ follow each other on $C^*$ in the forward direction. Since $|C_{A',j_k}|=2k+2$,  we find an even cycle $A'\overleftarrow{C^*}x_{j_2}A\overrightarrow{C^*}y_{j_k}A'$ of length $2k+4$, a contradiction. Hence, $e_0=2$, implying that $|C^*|=2k+3$.

(II) Suppose we have  $k$ even cycles $C_{A,j_1},\ldots,C_{A,j_k}$ or  $k$ even cycles $C_{A',j_1}',\ldots,C_{A',j_k}'$. By Claim \ref{A'e}, the position of $A$ and $A'$ is symmetry.  And we only need to consider the former case.
 Let $e_0=|x_{j_k}\overrightarrow{C^*}A|$. Then $e_0$ is even. 
Let $z=x_{j_1}^-$. Since $z\overleftarrow{C^*}Ax_{j_1}\overrightarrow{C^*}A_0PB$ is also a longest path in $G-V(C)$, by our choices (1) and (2) of $P$, $N(z)\ss C^*$.
Then there are $2k-1$ chords $z_jz$ in $C^*$ with $1\leq j\leq 2k-1$. 
 Without loss of generality, suppose $z_1,\ldots,z_{2k-1}$ follow each other on $C^*$ in the forward direction. 
 Let $H_z$ be the union of $C^*$ and $z_1z,\ldots,z_{2k-1}z$.
Applying Lemma \ref{C2k-1odd} on $H=H_z$,
 we have  
$k$ even cycles $C_{z,j}$ or $k$ even cycles $C_{z,j}'$.
If we have $k$ even cycles $C_{z,j_1},\ldots,C_{z,j_k}$ for some  $j_1<\cdots <j_k$ and $e_0\geq 4$, note that $|x_{j_k}\overrightarrow{C^*}A|=|z_{j_k}\overrightarrow{C^*}z|=e_0$,
then $A,z,z_{j_k}$ follow each other on $C^*$ in the forward direction.
We can find an even cycle $z\overleftarrow{C^*}Ax_{j_1}\overrightarrow{C^*}z_{j_k}z$ of length $2k+4$, a contradiction. 
If we have $k$ even cycles $C_{z,j_1}',\ldots,C_{z,j_k}'$ for some $j_1>\cdots >j_k$ and $e_0\geq 4$, then $x_{j_k},z_{j_2},z_{j_1},A$ follow each other on $C^*$ in the forward direction. We can  find an even cycle $A\overleftarrow{C^*}z_{j_{2}}z\overrightarrow{C^*}x_{j_k}A$ of length $2k+4$, a contradiction. Hence, $e_0=2$ in both cases, implying that $|C^*|=2k+3$.

By (I) and (II), $|C^*|=2k+3$. 
Since $d(A)\geq 2k+1$, denote the possible non-neighbor  of $A$ by $x$. Then $x\notin \{A^+,A,A_0\}$ and $x^-\notin \{A,A_0,A'\}$. For any vertex $v\in V(C^*)\backslash\{x^-,A,A_0,A'\}$, since $v\overleftarrow{C^*}Av^+\overrightarrow{C^*}A_0PB$ is also a longest path in $G-V(C)$, by our choices (1) and (2) of $P$, $N(v)\ss C^*$. Now we have $d_{C^*}(v)\geq 2k+1$ for any $v\in V(C^*)\backslash \{x^-,A_0\}$. Since $|V(C^*)|=2k+3$, every vertex $v\in V(C^*)\backslash \{x^-,A_0\}$ lies in a triangle $vv^+v^{2+}v$ or $vv^-v^{2-}v$. Note that $x^-A_0\notin E(C^*)$. Then for any $uv\in E(C^*)$, either $u$ or $v$ lies in a triangle 
 consisting of  three consecutive  vertices in $C^*$.

Since $G$ is 2-connected, 
we can find two disjoint paths $s_1P_1t_1$ and $s_2P_2t_2$ connecting $V(C)$ and $V(C^*)$, where $s_1,s_2\in V(C)$ and $t_1,t_2\in V(C^*)$.
Let $Q_1$ be a $(s_1,s_2)$-path in $C$ of length at least $k+1$.
Since $|V(C^*)|=2k+3$, we can find a   $(t_1,t_2)$-path $Q_2\ss C^*$  of length at least $k+2$. Since $|Q_2|\geq k+2\geq 5$, let $uv\ss Q_2$ and $d_{Q_2}(t_1,u)=d_{Q_2}(t_1,v)-1=2$. Since 
$u$ or $v$ lies in 
a triangle containing three consecutive  vertices in $C^*$, we can find a ($t_1,t_2$)-path $Q_2'$ in $G[V(Q_2)]$ of length $|Q_2|-1\geq k+1$.
Then we have an even cycle $s_1P_1t_1Q_2t_2P_2s_2Q_1s_1$ or $s_1P_1t_1Q_2't_2P_2s_2Q_1s_1$ of length larger than $2k+2$, a contradiction.

\vskip 3mm

Now we suppose that  there exists a longest path $P$ in $G-V(C)$ with two endpoints $A$ and $B$ such that $N(A)\cap C\not=\emptyset$ and $N(B)\cap C\not=\emptyset$. 
If possible, we choose such a longest path $P$  that $N(A)\cap C=N(B)\cap C$.
\vskip 1mm

\noindent{\bf{Case 4.}}
$N(A)\cap C\not=N(B)\cap C$. 
Let $P_x$ be any longest path in $G-V(C)$ with two endpoints $x,x'$ and both $N(x)\cap C$ and $N(x')\cap C$ are not empty. Then $N(x)\cap C\not=N(x')\cap C$. 
Let $|N(x)\cap V(C)|=q\geq 1$ and $|N(x)\cap V(P_x)|= p\geq 1$.  Since $d(x)\geq 2k+1$, $p+q\geq 2k+1$.  
Let $N(x)\cap V(C)=\{x_1,\ldots,x_q\}$ and $y\in N(x')\cap C\backslash N(x)$. Without loss of generality,  $y,x_1,\ldots,x_q$ follow each other on $C$
in the forward direction.
Let 
$N(x)\cap V(P_x)=\{z_{1},\ldots,z_{p}\}$ and $x,z_1,\ldots,z_p,x'$ follow each other on $P_x$ in this order.

Let $P_i=xz_iP_xx'$ and $Q_i=y\overrightarrow{C}x_i$. Then $P_x=P_1$. Let $Q_i'=y\overleftarrow{C}x_i$. Since $|C|=|Q_i'|+|Q_i|$, $|Q_i'|\equiv |Q_i| $.
Let $\mathcal{P}_{x}^1=\{P_i~|~|P_i|\equiv |P_1|\}$ and $\mathcal{P}_x^2=\{P_1,P_2,\ldots,P_p\}\backslash \mathcal{P}_x^1$. 
Let $\mathcal{Q}_{x}^1=\{Q_i~|~|Q_i|\equiv |P_1|\}$ and $\mathcal{Q}_x^2=\{Q_1,\ldots,Q_q\}\backslash \mathcal{Q}_{x}^1$. Then $|\mathcal{P}_{x}^1|+|\mathcal{P}_x^2|=p$ and $|\mathcal{Q}_{x}^1|+|\mathcal{Q}_x^2|=q$. 

\begin{claim}
    \label{Qi}
$|\mathcal{Q}_x^i|\leq k$ for any $i\in \{1,2\}$. 
\end{claim}

\noindent{\bf{Proof.}}
If $|\mathcal{Q}_x^i|\geq k+1$, denote $k+1$ of them by $Q_{j_0},\ldots,Q_{j_{k}}$, and $j_0<\cdots<j_{k}$. For any $0\leq i_1,i_2\leq k$, since 
$|Q_{j_{i_1}}|\equiv |Q_{j_{i_2}}|$, $|x_{j_{i_1}}\overrightarrow{C}x_{j_{i_2}}|$ is even. For any $i\leq k$, since  
we have an  even cycle $xx_{j_i}\overrightarrow{C}x_{j_{i'}}x$ for every $ i'\in \{0,\ldots,k\}-i$ and $|L_e(G)|=k$, we can check that
$|xx_{j_i}\overrightarrow{C}x_{j_{i-1}}x|=|C|$, implying that $|x_{j_{i-1}}\overrightarrow{C}x_{j_i}|=2$ for any $0\leq i\leq k$. Then $2|L_e(G)|+2=|C|=2k+2$. Since $y,x_1,\ldots,x_q$ follow each other on $C$
in the forward direction,  $x_{j_k}=x_q$, $x_{j_0}=x_1$ and $x_qyx_1\ss C$.
If $|P_x|$ is odd, then $xx_1\overrightarrow{C}yx'P_xx$ is an even cycle larger than $C$, a contradiction. 
Now we suppose $|P_x|$ is even. If $x_{j_i}\in N(x')$, then $xx_{j_{i+1}}\overrightarrow{C}x_{j_i}x'P_xx$ is an even cycle larger than $C$, a contradiction. If $x_{q}^-\in N(x')$ or $x_1^+\in N(x')$, then $xx_1\overrightarrow{C}x_{q}^-x'yx_{q}x$ or $xx_{j_1}\overrightarrow{C}yx'x_1^+x_{1}x$ is an even cycle of length $2k+4$, respectively, a contradiction. Then $|N(x')\cap C|\leq 2k+2-(k+3)=k-1$ and $|N(x')\cap P_x|\geq k+2$.
Denote $k+2$ of them by
$z_1',\ldots,z_{k+2}'$ and $x',z_1',\ldots,z_{k+2}',x$ follow each other on $P_x$ in this order.
If  $|x'z_i'P_xx|$ is odd, then  $xx_1\overrightarrow{C}yx'z_i'P_xx$ is an even cycle larger than $C$, a contradiction. Hence, $|x'z_i'P_xx|$ is even for every $1\leq i\leq k+2$. But we can find an even cycle $x'z_1'P_xz_i'x'$ for every $2\leq i\leq k+2$, a contradiction.
\bqed

\begin{claim}\label{P1}
$|\mathcal{P}_{x}^1|=k+1$ and $|\mathcal{Q}_{x}^1|=0$. As a corollary, $\mathcal{Q}_x^2=\{Q_1,\ldots,Q_q\}$, $|Q_i|\equiv |Q_i'|\not\equiv |P_x| $ and $|x_i\overrightarrow{C}x_j|$ is even  for any $1\leq i,j\leq q$.
\end{claim}

\noindent{\bf{Proof.}} Since $|\mathcal{P}_{x}^1|+|\mathcal{Q}_{x}^1|-1+|\mathcal{P}_x^2|+|\mathcal{Q}_x^2|-1=p+q-2\geq 2k-1$, suppose  $|\mathcal{P}_x^i|+|\mathcal{Q}_x^i|-1\geq k$ for some $i\in \{1,2\}$. 
Let $\mathcal{P}_x^i=\{P_{i_1},\ldots,P_{i_s}\}$ and $\mathcal{Q}_x^i=\{Q_{j_1},\ldots,Q_{j_t}\}$, where $i_1<\cdots<i_s$ and $j_1<\cdots<j_t$. Then $s+t\geq k+1$.
By Claim \ref{Qi}, $t\leq k$ and $s\geq 1$.

If $s=1$, then $t=k$. And  we have $k$ even cycles 
$xP_{i_1}x'yQ_{j_1}x_{j_1}x,\ldots,xP_{i_1}x'yQ_{j_k}x_{j_k}x$.
Then $|xP_{i_1}x'yQ_{j_k}x_{j_k}x|=|C|$ and 
$|Q'_{j_k}|=|P_{i_1}|+2$. 
We also have 
$ k$ even cycles 
$xP_{i_1}x'yQ_{j_k}'x_{j_k}x,\ldots,xP_{i_1}x'yQ_{j_1}'x_{j_1}x$. Then $|Q_{j_1}|=|P_{i_1}|+2$. And the smallest even cycle
$xP_{i_1}x'yQ_{j_1}x_{j_1}x$ has length $2|P_{i_1}|+4>4$.
Since these $k$ even cycles have consecutive even lengths, $|x_{j_1}\overrightarrow{C}x_{j_2}|=2$. But $xx_{j_1}\overrightarrow{C}x_{j_2}x$ is a 4-cycle smaller than $xP_{i_1}x'yQ_{j_1}x_{j_1}x$, a contradiction.

If $s\geq 2$ and $t\geq 1$, we have $t$ even cycles 
$xP_{i_s}x'yQ_{j_1}x_{j_1}x,\ldots,xP_{i_s}x'yQ_{j_t}x_{j_t}x$ and $s-1$ even cycles
$xP_{i_{s-1}}x'yQ_{j_t}x_{j_t}x,\ldots, xP_{i_{1}}x'yQ_{j_t}x_{j_t}x$.
Then we have $s+t-1= k$ even cycles and $|Q_{j_t}'|=|P_{i_1}|+2$.
We also have 
 even cycles 
$xP_{i_s}x'yQ_{j_t}'x_{j_t}x,\ldots,xP_{i_s}x'yQ_{j_1}'x_{j_1}x$ and 
$xP_{i_{s-1}}x'yQ_{j_1}'x_{j_1}x,\ldots, xP_{i_{1}}x'yQ_{j_1}'x_{j_1}x$, implying that $|Q_{j_1}|=|P_{i_1}|+2$. That means the lengths of all even cycles above are larger than $|P_{i_1}|+2$. But $xz_{i_1}P_xz_{i_2}x$ is an even cycle of length $|P_{i_1}|-|P_{i_2}|+2<|P_{i_1}|+2$, a contradiction to $|L_e(G)|= k$.

If $s\geq 2$ and $t=0$, then $s\geq k+1$.
We have  even cycles $xz_{i_1}P_xz_{i_2}x,\ldots,xz_{i_1}P_xz_{i_{s}}x$. Then $s=|\mathcal{P}_x^i|=k+1$. If $i=2$, then $P_x\notin\mathcal{P}_x^2$ and each of the above $k$ even cycles has length at most $|P_x|$. Since $t=|\mathcal{Q}_x^2|=0$, $\mathcal{Q}_x^1=\{Q_1,\ldots,Q_q\}$. But then $xP_xx'yQ_1x_{1}x$ is an even cycle larger than the above $k$ even cycles, a contradiction. Hence, $i=1$, that is, $s=|\mathcal{P}_x^1|=k+1$ and $t=|\mathcal{Q}_x^1|=0$.
 \bqed

\begin{claim}\label{z}
    $\mathcal{P}_x^2=\emptyset$. 
\end{claim}

\noindent{\bf{Proof.}}
By Claim \ref{P1}, let $z_{i_1},\ldots,z_{i_{k+1}}$ be $k+1$ neighbors of $x$ in $P_x$ such that $P_{i_j}\in \mathcal{P}_x^1$ for every $1\leq j\leq k+1$ and $i_1<\cdots<i_{k+1}$.  Then $i_1=1$, we have $k$ even cycles $xP_xz_{i_2}x,\ldots,xP_xz_{i_{k+1}}x$ and $|xP_xz_{i_{k+1}}x|=|C|$. Let $N(x')\cap V(C)=\{x_1',\ldots,x_{q'}'\}$, $N(x')\cap V(P)=\{z_{1}',\ldots,z_{p'}'\}$ and $y'\in N(x)\cap C\backslash N(x')$.
We can also define $\mathcal{P}_{x'}^i$ and $\mathcal{Q}_{x'}^i$ for $i\in \{1,2\}$.
Applying  Claims \ref{Qi} and \ref{P1} on $x=x'$, $x'$ also has $k+1$ neighbors $z'_{j_1},\ldots,z'_{j_{k+1}}$ in $P_x$ such that $j_1=1$,  we have $k$ even cycles $x'P_xz'_{j_2}x',\ldots,x'P_xz'_{j_{k+1}}x'$ of consecutive even lengths, $|x'P_xz'_{j_{k+1}}x'|=|C|$ and  $|x'P_xz'_{j_{k}}x'|=|C|-2$.

Suppose to the contrary that $\mathcal{P}_x^2\not=\emptyset$. By Claim \ref{P1},  there exists $z\in N(x)\cap V(P_x)$ such that $|xzP_xx'|\equiv |Q_1|\equiv |Q_1'|$. 
 Without loss of generality, suppose  $|Q_1|\geq \frac{|C|}{2}$. Since $xzP_xx'yQ_1x_1x$ is an even cycle, it must have length at most $|C|$. Then $|xzP_x{x'}|\leq \frac{|C|}{2}-2$. 
 Since $|x'P_xz'_{j_{k}}x'|=|C|-2$, $x,z'_{j_{k}},z,x'$ follow each other on $P_x$ in this order. Then we can check that $|xzP_xz_{j_{k}}'{x'}|=|xzP_x{x'}|+|x'P_xz'_{j_{k}}x'|-2|zP_x{x'}|\geq \frac{|C|}{2}+2$. Since    $|xzP_xz_{j_{k}}'{x'}|\equiv |xzP_x{x'}|\equiv |Q_1| $,
$xzP_xz_{j_{k}}'{x'}yQ_1x_1x$ is an even cycle of length at least $|C|+4$, a contradiction.\bqed

\vskip 1mm
Applying Claims \ref{Qi} - \ref{z} on $(P_x,x)=(P,A)$ or $(P,B)$, we have the following claim.

\begin{claim}\label{evenxij}

(1) $N(A)\cap P=\{z_1,\ldots,z_{k+1}\}$
and $N(B)\cap P=\{z_1',\ldots,z_{k+1}'\}$.

(2) $|Az_iPB|\equiv |APz_j'B|\equiv |P|$.

(3) $N(A)\cap C=\{x_1,\ldots,x_k\}$ and   
$N(B)\cap C=\{x_1',\ldots,x_{k}'\}$.

(4)
For any $1\leq i,j\leq k$, $|x_i\overrightarrow{C}x_j|$   and $|x_i'\overrightarrow{C}x_j'|$ are both even.  

(5) If $y\in N(B)\cap C\backslash N(A)$, then $|y\overrightarrow{C}x_i|\equiv |y\overleftarrow{C}x_i| \not\equiv |P| $.

(6) If $y'\in N(A)\cap C\backslash N(B)$, then $|y'\overrightarrow{C}x_i'|\equiv |y'\overleftarrow{C}x'_i| \not\equiv |P|$.

\end{claim}

Denote $z_2^-$ as the 
predecessor of $z_{2}$ in $APB$. Note that $P'=z_2^-PAz_{2}PB$ is also a longest path in $ G-V(C)$. If $N(z_2^-)\cap C=\emptyset$, then $N(z_2^-)\ss P$. Denote $2k+1$ neighbors of $z_2^-$ as $w_1,\ldots,w_{2k+1}$, and they follow each other on $P'=z_2^-PAz_{2}PB$ in this order. If $|z_2^-w_iP'B|\equiv |P'| $  for any $1\leq i\leq 2k+1$, then we have $2k$ even cycles $z_2^-P'w_2z_2^-,\ldots,z_2^-P'w_{2k+1}z_2^-$, a contradiction. If there exists $w\in \{w_1,\ldots,w_{2k+1}\}$
such that $|z_2^-wP'B|\not\equiv |P'| $, that is, $|z_2^-P'w|$ is even, we can check that either $A,w,z_2^-$ or $A,z_2^-,w$ folllow each other on $P$ in this order. In the former case, $|APwz_2^-PB|\not\equiv |P|$. In the later case, since $|z_2^-P'w|=|APz_2^-|+|z_2^-Pw|$ and both $|z_2^-P'w|$ and $|APz_2^-|$ 
are even, $|z_2^-Pw|$ is even. Then $|APz_2^-wPB|\not\equiv |P|$. In both cases, we have a $(A,B)$-path $P_z$ such that $|P_z|\not\equiv|P|$.
Let $y\in N(B)\cap C\backslash N(A)$, $Q_i=y\overrightarrow{C}x_i$ and $Q_i'=y\overleftarrow{C}x_i$.
By Claim \ref{evenxij}(5),
$|P_z|\equiv |Q_i| \equiv |Q_i'|$ for any $1\leq i\leq k$.
Now we have two sets of $k$ even cycles: $AP_zByQ_1x_1A,\ldots,AP_zByQ_kx_kA$ and $AP_zByQ_k'x_kA,\ldots,AP_zByQ_1'x_1A$.  Then $|C|=|AP_zByQ_kx_kA|=|AP_zByQ_1'x_1A|$,  $|Q_1|=|Q_k'|=|P_z|+2$, and the smallest even cycle has length $2|P_z|+4>4$. 
Since these even cycles have consecutive even lengths, $|x_1\overrightarrow{C}x_2|=2$. Then we get a 4-cycle $Ax_1\overrightarrow{C}x_2A$, a contradiction.

Now we suppose $N(z_2^-)\cap C \not=\emptyset$. Note that $P'=z_2^-PAz_{2}PB$ is also a longest path in $ G-V(C)$.  By our choice of $P$,
let $y'\in N(B)\cap C\backslash N(z_2^-)$.
Applying Claims \ref{Qi} - \ref{z} on $(x,x',y,P_x)=(z_2^-,B,y',P')$, $N(z_2^-)\cap C=\{w_1,\ldots,w_k\}$, $|w_i\overrightarrow{C} w_j|$ is even and $|y'\overrightarrow{C}w_i|\not\equiv |P|$ for any $1\leq i,j\leq k$.
By Claim \ref{evenxij}(4)(5), for any $1\leq i,j\leq k$,
$|x_i\overrightarrow{C} x_j|$ is even, $|y\overrightarrow{C}x_i|\not\equiv |P|$ and    $|y\overrightarrow{C}y'|$ is even. Then $|x_i\overrightarrow{C}w_j|$ is even for any $1\leq i,j\leq k$.

If $N(z_2^-)\cap C=N(A)\cap C$, note that $|APz_2^-|$ and  $|x_i\overrightarrow{C}x_j|$ is even for any $1\leq i,j\leq k$, we get $k$ even cycles $ APz_2^-x_1A,APz_2^-x_1\overrightarrow{C}x_2A,\ldots,APz_2^-x_1\overrightarrow{C}x_kA$ of consecutive even lengths, implying that $|x_1\overrightarrow{C}x_2|=2$.  But we have an even cycle $x_2\overrightarrow{C}x_1APz_2^-x_2$ larger than $C$, a contradiction.

If $N(z_2^-)\cap C\not=N(A)\cap C$,
without loss of generality, suppose  $w_1\in N(z_2^-)\cap C \backslash N(A)$ and $w_1,x_1,\ldots,x_k$ follow each other on $C$ in the forward direction. Since $|APz_2^-|$ and $|w_1\overrightarrow{C}x_i|$ is even, we have $k$ even cycles $APz_2^-w_1\overrightarrow{C}x_1A,\ldots,APz_2^-w_1\overrightarrow{C}x_kA$ of consecutive even lengths. Then  $|APz_2^-w_1\overrightarrow{C}x_1A|>4$ and  $|x_1\overrightarrow{C}x_2|=2$. But we have a 4-cycle $Ax_1\overrightarrow{C}x_2A$, a contradiction.

\vskip 3mm

\noindent{\bf{Case 5.}}
 $N(A)\cap C=N(B)\cap C$.
Let $x\in \{A,B\}$ and $x'\in \{A,B\}-x$.  
Let $N(x)\cap V(C)=\{x_1,\ldots,x_q\}$ and   they follow each other on $C$
in the forward direction.
Let 
$N(x)\cap V(P)=\{z_{1},\ldots,z_{p}\}$ and $x,z_1,\ldots,z_p,x'$ follow each other on $P$ in this order. Since $d(x)\geq 2k+1$, $p+q\geq 2k+1$.

\begin{claim}\label{q>1}
    $|N(x)\cap C|\geq 2$.
\end{claim}

\noindent{\bf{Proof.}} Suppose to the contrary that 
 $N(x)\cap C=\{x_0\}$.  
Let $C^*=x_0xPx'x_0$. Then   $C\cap C^*=\{x_0\}$ and  all  neighbors of both $x$ and $x'$ are in $C^*$. 
Since  $x'x_0x\ss C^*$ $N(x)\ss C^*$ and $N(x')\ss C^*$, the positions of $x$ and $x'$ are symmetric. 
Let $xx_j$ or $x'x_j'$ be a chord incident with $x$ or $x'$ for every $1\leq j\leq 2k-1$. Let $H_x$ be the union of $C^*$ and $xx_1,\ldots,xx_{2k-1}$ and $H_{x'}$ be the union of $C^*$ and $x'x'_1,\ldots,x'x'_{2k-1}$.
Since $|C^*|$ is even or odd, applying Lemma \ref{C2k-1} or \ref{C2k-1odd} on $H=H_x$, $|C^*|\geq |C|$. 
Define 
 $x'x_0xPx'$  of $C^*$ as the \emph{forward} direction, and denote by $\overrightarrow{C^*}$. And define another direction of $C^*$ as the \emph{backward} direction, and denote by $\overleftarrow{C^*}$.
Since $G$ is 2-connected, there is an ($u,v$)-path $P_{uv}$ connecting $C$ and $C^*$ in $G-x_0$, where $u\in C$ and $v\in C^*$. Without loss of generality, suppose $|u\overrightarrow{C}x_0|\geq \frac{|C|}{2}$ and $|x_0\overrightarrow{C^*}v|\geq \frac{|C|}{2}$. If $|u\overrightarrow{C}x_0|+|P_{uv}|\equiv |x_0\overrightarrow{C^*}v| $, then we find an even cycle $u\overrightarrow{C}x_0\overrightarrow{C^*}vP_{uv}u$ larger than $C$, a contradiction. Now we  consider that  $|u\overrightarrow{C}x_0|+|P_{uv}|\not\equiv |x_0\overrightarrow{C^*}v| $. In the rest of this proof, We aim to obtain a contradiction by finding another $(x_0,v)$-path in $G[V(C^*)]$ such that its length has a  different parity with $|x_0\overrightarrow{C^*}v|$ and is at least  
 $\frac{|C|}{2}$.
Note that $v\notin \{x,x_0,x'\}$.

If   $C^*$ is even, 
by applying  Lemma \ref{C2k-1} on $H=H_x$, we have $|C^*|=|C|=2k+2$ and  $|x\overrightarrow{C^*}x_1|=|x\overleftarrow{C^*}x_{2k-1}|=2$. Then $x_{2k-1}=x'$ and  $x_0,x,v,x_{2k-1}$ follow each other on $C^*$ in the forward direction.
Since $xx_0x_{2k-1}x$ is a triangle,  $|x_0x_{2k-1}x\overrightarrow{C^*}v|>|x_0\overrightarrow{C^*}v|$ and $|x_0x_{2k-1}x\overrightarrow{C^*}v|\not\equiv |x_0\overrightarrow{C^*}v|$. And $x_0x_{2k-1}x\overrightarrow{C^*}v$ is required.

If $C^*$ is odd, by applying  Lemma \ref{C2k-1odd} on $H=H_x$, then $|C^*|>|C|$ and we have  $k$ even cycle lengths $4,\ldots,2k+2$.
If there exists $x_j$ such that $x,v,x_j$ follow each other on $C^*$ in the forward direction and $C_{x,j}$ is even, then $C_{x,j}'$ is odd. And
$x_0\overleftarrow{C^*}x_jx\overrightarrow{C^*}v$ is a path as required, as shown in Figure \ref{x0v}(a).
Otherwise,  
let
$j_{k-2},j_{k-1}$ be two indices such  that $C_{x,j_{k-2}},C_{x,j_{k-1}}$  are two largest even cycles among all even $C_{x,j}$. Then $x,x_{j_{k-2}},x_{j_{k-1}},v$ or $x,x_{j_{k-2}},x_{j_{k-1}}=v$ follow each other on $C^*$ in the forward direction. 

Since we have $k$ even cycle lengths $4,\ldots,2k+2$, by Lemma \ref{C2k-1odd}, we can check   $|C_{x,j_{k-2}}|\in \{2k-2,2k\}$ and $|C_{x,j_{k-1}}|\in \{2k,2k+2\}$. Applying  Lemma \ref{C2k-1odd} on $H=H_{x'}$,
we  have at least $k-1$ even cycles $C_{x',j}$. Then there exists an index $j_a$ such that  $|C_{x',j_a}|\in \{|C_{x,j_{k-2}}|,|C_{x,j_{k-1}}|\}$ and $x',x,x'_{j_a},v$ follow each other on $C^*$ in the forward direction.
Then $|C_{x',j_a}|=2(k-i)+2$ for some $i\in \{0,1,2\}$. That is, $|x'\overrightarrow{C^*}x'_{j_a}|=2(k-i)+1$
and $|x\overrightarrow{C^*}x'_{j_a}|=2(k-i)-1$
for some $i\in \{0,1,2\}$.


If there exists $x_j$ such that $x,v,x_j$ follow each other on $C^*$ and $C_{x,j}$ is odd, then
$x_0x'x'_{j_a}\overleftarrow{C^*}xx_j\overleftarrow{C^*}v$, as shown in Figure \ref{x0v}(b), is an $(x_0,v)$-path of length $1+1+2(k-i)-1+1+|x_j\overleftarrow{C^*}v|\geq 2(k-2)+3=2k-1\geq \frac{|C|}{2}$.
Since $C_{x,j}$  is odd, $|x\overrightarrow{C^*}x_j|$ is even. Then $|x_0\overrightarrow{C^*}x_j|$ is odd, implying that $|x_j\overleftarrow{C^*}v|\not\equiv |x_0\overrightarrow{C^*}v|$. Since $|x_0x'x'_{j_a}\overleftarrow{C^*}xx_j\overleftarrow{C^*}v|\equiv |x_j\overleftarrow{C^*}v|$, we get $|x_0x'x'_{j_a}\overleftarrow{C^*}xx_j\overleftarrow{C^*}v|\not\equiv |x_0\overrightarrow{C^*}v|$. And $x_0x'x'_{j_a}\overleftarrow{C^*}xx_j\overleftarrow{C^*}v$ is  required.

\begin{figure}[!htb]
\centering
\includegraphics[height=0.3\textwidth]{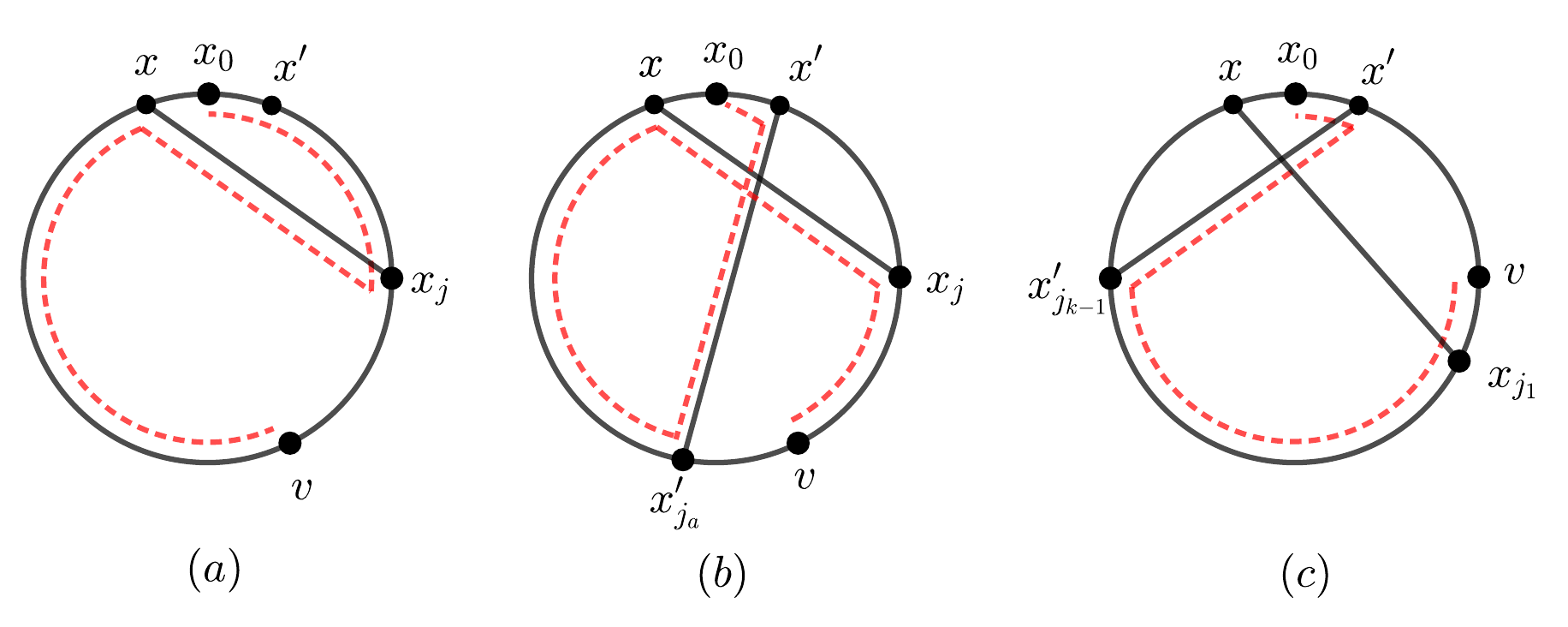}
\caption{Some $(x_0,v)$-paths}
\label{x0v}
\end{figure}

Now we suppose all $x_j$ lie in $x\overrightarrow{C^*}v$.
Let ${j_1}$ be an index such that $C_{x,j_1}'$
is the smallest even cycle among all $C_{x,j}'$. Then $x_{j_1}$ lies in  $x\overrightarrow{C^*}v$
 and $|x\overleftarrow{C^*}x_{j_1}|\in \{3,5\}$, implying that $|x\overleftarrow{C^*}v|\leq 5$ and  $|x'\overleftarrow{C^*}v|\leq 3$. Let $j_{k-1}$ be an index such that $C_{x',j_{k-1}}'$
is the largest even cycle among all $C_{x',j}'$. Then $x',x'_{j_{k-1}},v$ follow each other on $C^*$ in the forward direction and 
$|C_{x',j_{k-1}}'|=2(k-i)+2$ for some $i\in \{0,1\}$. Since $|x'\overleftarrow{C^*}x'_{j_{k-1}}|= 2(k-i)+1$ for some $i\in \{0,1\}$, we get  $|x'_{j_{k-1}}\overrightarrow{C^*}v|\geq 2(k-i)+1-3\geq 2k-4$.
Then $x_0x'x'_{j_{k-1}}\overrightarrow{C^*}v$ is a path of length $2+2k-4=2k-2\geq \frac{|C|}{2}$.
Since $C_{x',j_{k-1}}$ is an odd cycle, we get $|x_0\overrightarrow{C^*}v|\not\equiv |x_0x'x'_{j_{k-1}}\overrightarrow{C^*}v|$. And $x_0x'x'_{j_{k-1}}\overrightarrow{C^*}v$ is required, as shown in Figure \ref{x0v}(c).
\bqed

\vskip 2mm
By Claim \ref{q>1},  $|N(x)\cap C|=q\geq 2$. Then we have a 4-cycle $xx_1x'x_2x$. Since we have $k$ consecutive even cycle lengths and $|L_e(G)|=k$, we have even cycle lengths $4,6,\ldots,2k+2$. And $|C|=2k+2$.  Let $P_i=xz_iPx'$ and $Q_i=x_1\overrightarrow{C}x_i$. Then $P_1=P$ and $|Q_1|=0$.  For $i\geq 2$, let $Q_i'=x_1\overleftarrow{C}x_i$. Since $|C|=|Q_i'|+|Q_i|$,  $|Q_i'|\equiv |Q_i| $.
Let $\mathcal{P}_{x}^1=\{P_i~|~|P_i|\equiv |P_1|\}$ and $\mathcal{P}_x^2=\{P_1,P_2,\ldots,P_p\}\backslash \mathcal{P}_x^1$. 
Let $\mathcal{Q}_{x}^1=\{Q_i~|~|Q_i|\equiv |P_1|\}$ and $\mathcal{Q}_x^2=\{Q_1,\ldots,Q_q\}\backslash \mathcal{Q}_{x}^1$. Then $|\mathcal{P}_{x}^1|+|\mathcal{P}_x^2|=p$ and $|\mathcal{Q}_{x}^1|+|\mathcal{Q}_x^2|=q$. Since $|\mathcal{P}_{x}^1|+|\mathcal{Q}_{x}^1|-1+|\mathcal{P}_x^2|+|\mathcal{Q}_x^2|-1=p+q-2\geq 2k-1$, one of $\{|\mathcal{P}_x^1|+|\mathcal{Q}_{x}^1|-1,|\mathcal{P}_x^2|+|\mathcal{Q}_x^2|-1\}$ is at least $k$.
Suppose $|\mathcal{P}_x^i|+|\mathcal{Q}_{x}^i|-1\geq k$ for some $i\in \{1,2\}$.
Let $\mathcal{P}_x^i=\{P_{i_1},\ldots,P_{i_s}\}$ and $\mathcal{Q}_{x}^i=\{Q_{j_1},\ldots,Q_{j_t}\}$. Then $s+t\geq k+1$.
\vskip 1mm

\noindent{\bf{Subcase 5.1.}} $s\geq 1$ and $t\geq 1 $ for some $x\in \{A,B\}$.
Note that $xP_{i_s}x'x_1Q_{j_1}x_{j_1}x,\ldots, $  $xP_{i_s}x'x_1Q_{j_t}x_{j_t}x$ and  $xP_{i_{s-1}}x'x_1Q_{j_t}x_{j_t}x,\ldots,xP_{i_1}x'x_1Q_{j_t}x_{j_t}x$ are even cycles. Then we have $s+t-1=k$ even cycles of lengths $4,6,\ldots,2k+2$.
If $|Q_{j_1}|\not=0$,  then $|P_{i_s}|=|Q_{j_1}|=1$ and $x_{j_1}=x_2$. But $x_1xP_{i_s}x'x_2\overrightarrow{C}x_1$ is an even cycle of length $2k+4$, a contradiction.
If $t\geq 2$ and $|Q_{j_1}|=0$, then $x_1=x_{j_1}$, every path in $\mathcal{P}_x^i$ has  even lengths and $|x_{1}\overrightarrow{C}x_{j_2}|=2$.
But  $x_1xP_{i_s}x'x_{j_2}\overrightarrow{C}x_1$ is an even cycle larger than $C$, a contradiction.

Now we suppose $t=1$ and $|Q_{j_1}|=0$. That is, $x_{j_1}=x_1$, $s\geq k$ and every path in $\mathcal{P}_x^i$ has  even lengths. Then we have $s=k$ even cycles $xP_{i_s}x'x_1x,\ldots,xP_{i_1}x'x_1x$ of lengths $4,\ldots,2k+2$.
Since $t=1$, all paths in $\{Q_2,\ldots,Q_q\}$ are odd.
If $|N(x)\cap C|\geq 3$, since $|Q_2|\equiv |Q_3|$, $|x_2\overrightarrow{C}x_3|$ is even. But $xP_{i_1}x'x_2\overrightarrow{C}x_3x$ is an even cycle larger than $xP_{i_1}x'x_1x$, a contradiction. If $|N(x)\cap C|=2$,  note that $|x_1\overrightarrow{C}x_2|\equiv|x_1\overleftarrow{C}x_2|\equiv1$, without loss of generality, suppose $|x_1\overrightarrow{C}x_2|\geq k+1$. Since $|N(x)\cap P|\geq 2k-1$, we have at least $k-1$ neighbors of $x$, $z_{i_1'},\ldots,z_{i_{k-1}'}$, such that $xz_{i_1'}Px',\ldots,xz_{i_{k-1}'}Px'$ have different odd lengths. 
Then the longest odd path, denote as $P'$, must have length at least $2k-3$. But we have an even cycle $xP'x'x_1\overrightarrow{C}x_2x$ larger than $C$, a contradiction.

\vskip 2mm

\noindent{\bf{Subcase 5.2.}}
$s=0$ for some $x\in \{A,B\}$. Since $P\in \mathcal{P}_x^1$, $s=|\mathcal{P}_x^2|=0$ and $t=|\mathcal{Q}_x^2|\geq k+1$.
 Note that the lengths of all paths in $\mathcal{Q}_x^2$ have the same parity. 
Then we have $t-1=k$ even cycles $xx_{j_1}\overrightarrow{C}x_{j_{2}}x,\ldots,xx_{j_1}\overrightarrow{C}x_{j_{k+1}}x$ of  lengths
$4,\ldots,2k+2$. And $|x_{j_i}\overrightarrow{C}x_{j_{i+1}}|=2$.
 If  $|P|$ is even, then $x_{j_1}xPx'x_{j_2}\overrightarrow{C}x_{j_1}$ is an even cycle larger than $C$, a contradiction.

Now we suppose that all paths in $\mathcal{P}_x^1$ are odd. If there exists a vertex $z\in N(x)\cap C\backslash\{x_{j_1},\ldots,x_{j_{k+1}}\}$, then $x_{j_i}zx_{j_{i+1}}\ss C$ for some $i$. But $x_{j_i}xPx'z\overrightarrow{C}x_{j_i}$ is an even cycle larger than $C$, a contradiction.
Then $|N(x)\cap P|\geq k\geq 3$. Let $z_2^-$ be the predecessor of $z_2$ in $xPx'$. Then $P'=z_2^-Pxz_2Px'$ is also a longest path in $G-V(C)$. 
Suppose $N(z_2^-)\cap C\not=\emptyset$.
Note that $|xPz_2^-|$ is even and $|z_2^-Px'|$ is odd.
If there exists a vertex $z\in N(z_2^-)\cap \{x_{j_1},\ldots,x_{j_{k+1}}\}$, let $z=x_{j_i}$ for some $i$. Then $x_{j_i}z_2^-Pxx_{j_{i+1}}\overrightarrow{C}x_{j_i}$ is an even cycle larger than $C$, a contradiction. If  there exists a vertex $z\in N(z_2^-)\cap C$ and $x_{j_i}zx_{j_{i+1}}\ss C$ for some $i$, then 
$x_{j_i}x'Pz_2^-z\overrightarrow{C}x_{j_i}$ is an even cycle larger than $C$, a contradiction.

If $N(z_2^-)\ss P$, denote $2k+1$ of them by $z_1',\ldots,z_{2k+1}'$, and they follow each other on $z_2^-P'x'$ in this order.
If all $z_2^-z_i'P'x'$ have odd lengths, then we have $2k$ even cycles $z_2^-P'z_2'z_2^-,\ldots,z_2^-P'z_{2k+1}'z_2^-$, a contradiction.
Otherwise, let $z$ be a neighbor of $z_2^-$ and $|z_2^-P'z|$ is even. 
We can check that either $x,z,z_2^-$ or $x,z_2^-,z$ follow each other on $P$ in this order. In the former case, $xPzz_2^-Px'$ has even length. In the later case, since $|z_2^-P'z|=|xPz_2^-|+|z_2^-Pz|$ and both $|z_2^-P'z|$ and $|xPz_2^-|$ 
are even, $|z_2^-Pz|$ is even. And $xPz_2^-zPx'$ has even length. In both cases, we have a $(x,x')$-path $P_z$ of even length. But then $x_{j_1}xP_zx'x_{j_2}\overrightarrow{C}x_{j_1}$ is an even cycle larger than $C$, a contradiction.

\vskip 2mm

\noindent{\bf{Subcase 5.3.}} $t=0$ for any $x\in \{A,B\}$. Then $|\mathcal{Q}_x^i|=0$ and $s=|\mathcal{P}_x^i|\geq k+1$. Since $|Q_1|=0$,  $Q_1\in \mathcal{Q}_x^{3-i}$. Then all paths in $\mathcal{Q}_x^{3-i}$ must have even lengths and all paths in $\mathcal{P}_x^i$ must have odd lengths.
Since we have $s-1$ even cycles $xz_{i_1}Pz_{i_2}x,\ldots,xz_{i_1}Pz_{i_s}x$, $s-1=k$ and all  have lengths at most $|P|+1$. If $|P|$ is even, then $xPx'x_1x$ is an even cycle of length $|P|+2$, a contradiction. Hence, $|P|$ is odd, $P=P_{i_1}$ and $ \mathcal{P}_x^i=\mathcal{P}_x^1=\{P_{i_1},\ldots,P_{i_{k+1}}\}$. Then $|\mathcal{Q}_x^1|=0$,  $\mathcal{Q}_x^2=\{Q_1,\ldots,Q_q\}$,  and $|x_i\overrightarrow{C}x_j|$ is even for any $1\leq i,j\leq q$.
Note that we have $k$ even cycles $xPz_{i_{2}}x,\ldots,xPz_{i_{k+1}}x$ of lengths $4,\ldots,2k+2$ and  $|C|=2k+2$.

Let $C^*=Bx_1APB$.  Then $C^*$ is an odd cycle larger than $C$ and $C^*\cap C=\{x_1\}$.
Define 
 $Bx_1APB$  of $C^*$ as the \emph{forward} direction, and denote by $\overrightarrow{C^*}$. And define another direction of $C^*$ as the \emph{backward} direction, and denote by $\overleftarrow{C^*}$.
Without loss of generality, suppose $|x_1\overrightarrow{C}x_2|\geq k+1$. 
Since $|x_1\overrightarrow{C}x_2|$ is even, $|Ax_1\overrightarrow{C}x_2B|\geq k+3$ and  is even. Taking $x=A$, then $A$ has $k+1$ neighbors in $P$ such that $AzPB$ has  odd length for every $z$ in these $k+1$ vertices.

If $A$ has a neighbor $z$ in $P$ such that $AzPB$ is even and $|AzPB|\geq k+1$, then $Ax_1\overrightarrow{C}x_2BPzA$ is an even cycle larger than $C$, a contradiction.
Suppose $A$ has a neighbor $z$ in $P$ such that $AzPB$ is even and $|AzPB|\leq k$. By taking $x=B$, there are $k+1$ neighbors of $B$ in $P$, denoted as $z_{i_1},\ldots,z_{i_{k+1}}$ such that $BPz_{i_{2}}B,\ldots,BPz_{i_{k+1}}B$ have lengths $4,\ldots,2k+2$.
Then  $|BPz_{i_{k}}|=2k-1$. Since 
$|zPB|\leq k-1$ and is odd, $A,z_{i_k},z,B$ follow each other on $P$ in this order. We can check that $Bz_{i_{k}}PzA$ is a path of even length $1+|z_{i_{k}}PB|-|zPB|+1\geq k+2$. Then $Ax_1\overrightarrow{C}x_2Bz_{i_{k}}PzA$ is an even cycle larger than $C$, a contradiction.

Now we remain to consider that $N(A)\cap P=\{z_1,\ldots,z_{k+1}\}$, $|Az_iPB|$ is odd and $N(A)\cap C=\{x_1,\ldots,x_k\}$.
Since $|x_i\overrightarrow{C}x_{i+1}|$ is even and $|C|=2k+2$, without loss of generality, let $|x_1\overrightarrow{C}x_{2}|=\cdots=|x_{k-1}\overrightarrow{C}x_{k}|=2$ and $|x_k\overrightarrow{C}x_1|=4$.
Let $z\in N(A)\cap P-z_1$ and $z^-$ be the predecessor of $z$ in $APB$. Then $z^-PAzPB$ is also a longest path in $G-V(C)$. Denote $x_{k+1}$ as the vertex in $C$ and $|x_k\overrightarrow{C}x_{k+1}|=|x_{k+1}\overrightarrow{C}x_1|=2$.

If $N(z^-)\cap C\not= \emptyset$, let $z'\in N(z^-)\cap C $. Note that $|APz^-|$ is even and $|z^-PB|$ is odd.  
If $z'=x_i$ for some $i$, then $x_iz^-PAx_{i+1}\overrightarrow{C}x_i$ or $x_iz^-PAx_{i-1}\overleftarrow{C}x_i$ is an even cycle larger than $C$, a contradiction. If $z'\notin \{x_1,\ldots,x_k,x_{k+1}\}$, then $z'x_i\ss C$ for some $i\leq k$. And either $z'z^-PBx_i\overrightarrow{C}z'$ or $z'z^-PBx_i\overleftarrow{C}z'$ is an even cycle larger than $C$, a contradiction.

If $N(z^-)\ss P$, denote $2k+1$ of them by $z_1',\ldots,z_{2k+1}'$, and they follow each other on $z^-P'x'$ in this order.
If all $z^-z_i'P'x'$ have odd lengths, then we have $2k$ even cycles $z^-P'z_2'z^-,\ldots,z^-P'z_{2k+1}'z^-$, a contradiction.
Otherwise, let $z$ be a neighbor of $z^-$ and $|z^-P'z|$ is even. 
We can check that either $x,z,z^-$ or $x,z^-,z$ follow each other on $P$ in this order. In the former case, $xPzz^-Px'$ has even length. In the later case, since $|z^-P'z|=|xPz^-|+|z^-Pz|$ and both $|z^-P'z|$ and $|xPz^-|$ 
are even, $|z^-Pz|$ is even. And $xPz^-zPx'$ has even length. In both cases, we have a $(x,x')$-path $P_z$ of even length. But then $x_{1}xP_zx'x_{2}\overrightarrow{C}x_{1}$ is an even cycle larger than $C$, a contradiction.  \bqed

\section{Some lemmas}

\begin{lem}\label{2-}
    $\chi(K_n^{2-})\leq n-2$.
\end{lem}

\noindent{\bf{Proof.}}
Let $uv,xy$ be two edges deleted in $K_n$. Then $W=K_n^{2-}[V(K_n^{2-})\backslash\{u,x\}]$ is a clique  of size $n-2$. Assign each vertex in $W$  the color $1,\ldots,n-2$, respectively. And assign $u$ the same color
of $v$, $x$ the same color
of $y$. Then this is a proper coloring and the conclusion holds.
\bqed

\begin{lem}\label{C'} 
Let $G$ be a 2-connected graph.   Let $C^*$ be an odd cycle  and $C^*$ has  at most one common vertex with a longest odd cycle $C$. Then $|C^*|<|C|$. 
\end{lem}

\noindent{\bf{Proof.}} The case when $C\cap C^*=\emptyset $ has been proved by Gy\'arf\'as \cite{G}. Now we suppose that $C^*\cap C=\{v\}$ and $|C^*|=|C|$. Since $G$ is 2-connected, there is a ($x,x^*$)-path $P$ connecting $C-v$ and $C^*-v$, where $x\in C-v$ and $x^*\in C^*-v$.
Since $C$ and $C^*$ are odd, let $Q_1,Q_2$ be two $(x,v)$-paths in $C$ and $Q_1^*,Q_2^*$ be two $(y,v)$-paths in $C^*$ such that $|Q_1|\equiv|Q_1^*|\equiv 1$ and $|Q_2|\equiv|Q_2^*|\equiv 0$.  If $|P|$ is even, then $xQ_1vQ_2^*x^*Px$ and $xQ_2vQ_1^*x^*Px$ are odd cycles. If $|P|$ is odd, then $xQ_1vQ_1^*x^*Px$ and $xQ_2vQ_2^*x^*Px$ are odd cycles. Since $|C|+|C^*|+2|P|\geq 2|C|+2$, we can check that one of $\{xQ_1vQ_1^*x^*Px,xQ_2vQ_2^*x^*Px\}$ or $\{xQ_1vQ_2^*x^*Px,xQ_2vQ_1^*x^*Px\}$ is odd but larger than $C$, a contradiction.\bqed

\begin{lem}\label{non-bipartite}
    Let $H$ be a   graph obtained by adding $2t-2$ chords in a cycle $C_0$ incident with the same vertex $x$. If $H$ is non-bipartite, then $|L_o(H)|\geq t$.  
\end{lem}

\noindent{\bf{Proof.}} 
We define one direction of $C_0$ as the \emph{forward} direction, and denote as $\overrightarrow{C_0}$. And define another direction of $C_0$ as the \emph{backward} direction, and denote as $\overleftarrow{C_0}$.

We prove it by induction on $t$. If $t=1$, then $H=C_0$ is an odd cycle, implying that $|L_o(H)|\geq 1$. Suppose that this lemma is true for $t'<t$, $t\geq 2$. 
Denote $2t-2$ chords as $xx_j$ with $j=1,\ldots,2t-2$, and $x_1,\ldots,x_{2t-2}$ follow each other on $C_0$ in the forward direction. 
If $C_0$ is an odd cycle,  each chord $xx_{j}$ separates the odd cycle $C_0$ into two small cycles $C_{1,j}=x\overrightarrow{C_0}x_{j}x$ and $C_{2,j}=x\overleftarrow{C_0}x_{j}x$. One is an even cycle and the other  is an odd cycle. Since $C_0$ has $2t-2$ chords $xx_j$, there are at least $t-1$ odd cycles $C_{1,j}$ of different lengths, or at least $t-1$ odd cycles $C_{2,j}$ of different lengths. Together with $C_0$, we get $|L_o(H)|\geq t$. 

If $C_0$ is an even cycle, then
each chord $xx_{j}$ separates the even cycle $C_0$ into two small cycles $C_{1,j}=x\overrightarrow{C_0}x_{j}x$ and $C_{2,j}=x\overleftarrow{C_0}x_{j}x$ of the same parity length. We say that $x_j$ is \emph{even} if both $C_{1,j}$ and $C_{2,j}$ are even cycles, and  $x_j$ is \emph{odd} if both $C_{1,j}$ and $C_{2,j}$ are odd cycles. If there exists $t$ odd $x_j$, then $|L_o(H)|\geq t$.

If the number of odd $x_j$ is at most $t-1$, then we have at least $t-1$ even $x_j$.
Since $H$ is non-bipartite, there must exist odd $x_j$ for some $j\leq 2t-2$.
Let $p$ be the smallest integer such that $x_p$ or $x_{2t-1-p}$ is odd, and $p<2t-1-p$. Then $p\leq t-1$ and $x_1,\ldots,x_{p-1},x_{2t-p},\ldots,x_{2t-2}$ are even. If $p=1$, without loss of generality, we suppose $x_1$ is odd.
 Since we have at least $t-1$ even $x_j$, we have at least $t-1$ odd cycles $xx_1\overrightarrow{C_0}x_jx$. Together with $xx_1\overrightarrow{C_0}x$, we get $|L_o(H)|\geq t$.

If $2\leq p\leq t-1$, let $H^*:=H[ V(xx_p\overrightarrow{C_0}x_{2t-1-p}x)]$.
If $H^*$ is non-bipartite, note that $x$ has $2(t-p)-2$ chords in $H^*$, by the induction hypothesis, then $|L_o(H^*)|\geq t-p$. Without loss of generality, suppose that $x_p$ is odd. Since $x_j$ is even for  every $j\in \{2t-p,\ldots,2t-2\}$, $xx_p\overrightarrow{C_0}x_jx$ is an odd cycle. Note that $xx_p\overrightarrow{C_0}x$ is also an odd cycle.
Then $|L_o(H)|\geq t-p+p-1+1=t$.

If $H^*$ is a bipartite graph, then $|x_p\overrightarrow{C_0}x_{2t-1-p}|$ is even. Since $C_0$ is an even cycle, by the definitions of $x_p$ and $x_{2t-1-p}$, $|x\overrightarrow{C_0}x_{p}|$ and $|x\overleftarrow{C_0}x_{2t-1-p}|$ are both even,  implying that
$x_p$ and $x_{2t-1-p}$ are odd.
Indeed,  $x_j$ is odd for every $p\leq j\leq 2t-1-p$.
Since $p\geq 2$, let $a_1=|x_{p-1}\overrightarrow{C_0}x_p|$ and 
$a_{2}=|x_{2t-1-p}\overrightarrow{C_0}x_{2t-p}|$.
Without loss of generality, suppose that $a_1\geq a_2$.
Since $x_{2t-p}$ is even, 
$C_j=xx_j\overrightarrow{C_0}x_{2t-p}x$ is an odd cycle for every $j\in \{p,\ldots,2t-1-p\}$.
For every  $j\leq p-1$, $x_j$ is even. Then $C_j=xx_j\overrightarrow{C_0}x_{2t-1-p}x$ is an odd cycle for $j\leq p-1$. Since $a_1\geq a_2$,
 $|C_{2t-1-p}|<\cdots<|C_p|\leq |C_{p-1}|<\cdots <|C_1|$, and
 we have found $2t-p-2$ odd cycles of different lengths. 
Together with $x\overrightarrow{C_0}x_{2t-1-p}x$, we find $2t-p-1\geq t$ odd cycles of different lengths, completing the proof. \bqed

\begin{lem}
    
\label{bipartite}
    Let $C'$ be an even cycle with $xx_0x'\ss C'$. 
    Let
    $H$ be a graph obtained by adding chords $xx_1,\ldots,xx_{2k-2}$   and  $x'x'_1,\ldots,x'x'_{2k-2}$ in $C'$.  If $k\geq 3$ and $H$ is bipartite, then for any two vertices $y,z\in H$, we can find $k+1$ $(y,z)$-paths of different  but  the same parity lengths in $H$. 
\end{lem}

\noindent{\bf{Proof.}}
  Since $H$ is bipartite,
all $(y,z)$-paths  have the same parity lengths, and we only need to find $k+1$ $(y,z)$-paths of different lengths.
 In most situations, we will find them in  the union of $C'$ and $xx_1,\ldots,xx_{2k-2}$.   
 In one special situation in the end of this proof, we will find them in  $H$. We define one direction of $C'$ as the \emph{forward} direction, and denote as $\overrightarrow{C'}$. And define another direction of $C'$ as the \emph{backward} direction, and denote as $\overleftarrow{C'}$. 
 Let $x_1,\ldots,x_{2k-2}$ follow each other on $C'$ in the forward direction.

If $\{x,x'\}= \{y,z\}$, without loss of generality, suppose $x=y$,  $x'=z$, and $x_0,y,x_1,$  $\ldots,x_{2k-2},z$ follow each other on $C'$ in the forward direction. 
Then there exist $2k-2\geq k+1$ $(y,z)$-paths $yx_1\overrightarrow{C'}z,\ldots, yx_{2k-2}\overrightarrow{C'}z$.

If one of $\{x,x'\}$ does not belong to $ \{y,z\}$,  
without loss of generality, suppose $x\notin \{y,z\}$ and  $x,y,z$ follow each other on $C'$ in the forward direction. We first consider that $x,y,x_1$ follow each other on $C'$ in the forward direction.
If $x,y,z,x_1$ or $x,y,z=x_1$ follow each other on $C'$ in the forward direction, then we find $2k-2\geq k+1$ $(y,z)$-paths $y\overleftarrow{C'}xx_1\overleftarrow{C'}z,\ldots,y\overleftarrow{C'}xx_{2k-2}\overleftarrow{C'}z$ of different lengths. If $x,y,x_1,z$ follow each other on $C'$ in the forward direction, let $x_1,\ldots, x_t$ be the internal vertices in $y\overrightarrow{C'}z$. Then we find $2k-2\geq k+1$ $(y,z)$-paths $y\overrightarrow{C'}x_1xx_{t+1}\overleftarrow{C'}z, \ldots,y\overrightarrow{C'}x_1xx_{2k-2}\overleftarrow{C'}z, y\overrightarrow{C'}x_1x\overleftarrow{C'}z,\ldots,y\overrightarrow{C'}x_tx\overleftarrow{C'}z$ of different lengths.

Now we suppose that $x,x_1,y$ or $x,x_1=y$ follow each other on $C'$ in the forward direction. If all $x_j$ lie in $x\overrightarrow{C'}y$, then we find $2k-2\geq k+1$ $(y,z)$-paths $y\overleftarrow{C'}x_{2k-2}x\overleftarrow{C'}z,\ldots,y\overleftarrow{C'}x_1x\overleftarrow{C'}z$  of different lengths.
If there are some $x_j$ not in $x\overrightarrow{C'}y$, 
suppose  $x_1,\ldots,x_s\in x\overrightarrow{C'}y$ and $1\leq s<2k-2$. If $x_{s+1},\ldots,x_{2k-2}$ are the internal vertices in 
$y\overrightarrow{C'}z$, then we can find   $2k-2\geq k+1$ $(y,z)$-paths $y\overleftarrow{C'}x_sxx_{2k-2}\overrightarrow{C'}z,\ldots,y\overleftarrow{C'}x_1xx_{2k-2}\overrightarrow{C'}z$, $y\overleftarrow{C'}xx_{2k-2}\overrightarrow{C'}z,\ldots,y\overleftarrow{C'}xx_{s+1}\overrightarrow{C'}z$   of different lengths.
If $x_{s+1},\ldots,x_{2k-2}$ lie in 
$z\overrightarrow{C'}x$, then we can find  $2k-2\geq k+1$ $(y,z)$-paths  $y\overleftarrow{C'}x_sxx_{s+1}\overleftarrow{C'}z,\ldots,y\overleftarrow{C'}x_1xx_{s+1}\overleftarrow{C'}z$, $y\overleftarrow{C'}xx_{s+1}\overleftarrow{C'}z,\ldots,y\overleftarrow{C'}xx_{2k-2}\overleftarrow{C'}z$ of different lengths.

Now we suppose that $x_{s+1},\ldots,x_t$ are the internal vertices in $y\overrightarrow{C'}z$ with $1\leq s<t<2k-2$.
We will find two sets $S_1,S_2$ of $(y,z)$-paths, and all the paths in each set  have different lengths.
The first set $S_1$ contains $t$ $(y,z)$-paths $y\overleftarrow{C'}x_sxx_t\overrightarrow{C'}z,\ldots, y\overleftarrow{C'}x_1xx_t\overrightarrow{C'}z,$ $ y\overleftarrow{C'}xx_t\overrightarrow{C'}z,\ldots,y\overleftarrow{C'}xx_{s+1}\overrightarrow{C'}z$. The second set $S_2$ contains $s+2k-2-t+1=2k-t+s-1$  $(y,z)$-paths  $y\overleftarrow{C'}x_sxx_{t+1}\overleftarrow{C'}z,\ldots,y\overleftarrow{C'}x_1xx_{t+1}\overleftarrow{C'}z,y\overleftarrow{C'}xx_{t+1}\overleftarrow{C'}z,\ldots,$  $ y\overleftarrow{C'}xx_{2k-2}\overleftarrow{C'}z,y\overleftarrow{C'}z$. Note that $|S_1|+|S_2|=2k-1+s$ and $s\geq 1$. If $s\geq 2$, or $s=1$ and $t\not=k$, then we can check that one of these two sets has at least $k+1$ paths of different lengths. 
If $s=1$ and $t=k\geq 4$, then we have $2k-3\geq k+1$ $(y,z)$-paths 
$y\overrightarrow{C}x_{2}xx_{k+1}\overleftarrow{C'}z\ldots,y\overrightarrow{C}x_{2}xx_{2k-2}\overleftarrow{C'}z,y\overrightarrow{C'}x_2x\overleftarrow{C'}z, \ldots, y\overrightarrow{C'}x_{k}x\overleftarrow{C'}z$   of different lengths.

Now we remain to consider that $s=1$ and $t=k=3$ for $x$, as shown in Figure \ref{k+1_paths}.  
Let $O_1=|x\overrightarrow{C'}x_1|$, 
$O_2=|x_4\overrightarrow{C'}x|$,
$p_1=|x_1\overrightarrow{C'}y|$,
$p_2=|y\overrightarrow{C'}x_2|$,
$q_1=|x_3\overrightarrow{C'}z|$,
$q_2=|z\overrightarrow{C'}x_4|$ and $e_0=|x_2\overrightarrow{C'}x_3|$.

\begin{figure}[!htb]
\centering
\includegraphics[height=0.3\textwidth]{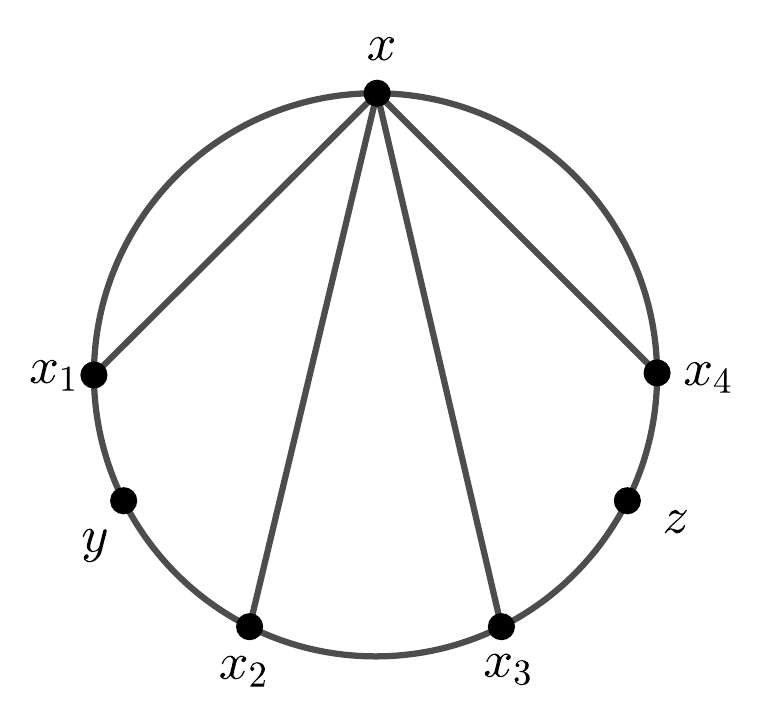}
\caption{$s=1$ and $t=k=3$.}
\label{k+1_paths}
\end{figure}

We first construct three $(y,z)$-path sets with their path length sets:
  
 $\{y\overleftarrow{C'}z,    y\overleftarrow{C'}x_1x\overleftarrow{C'}z,   y\overleftarrow{C'}x_1xx_4\overleftarrow{C'}z\}$ with 
 
\hspace{2mm} $\{p_1+O_1+O_2+q_2, p_1+1+O_2+q_2,p_1+2+q_2\}$, 
 \vskip 1mm

 $\{y\overleftarrow{C'}xx_2\overrightarrow{C'}z,    y\overleftarrow{C'}x_1xx_2\overrightarrow{C'}z,   y\overleftarrow{C'}x_1xx_3\overrightarrow{C'}z\}$  with 
 
\hspace{2mm}  $\{p_1+O_1+1+e_0+q_1, p_1+2+e_0+q_1,p_1+2+q_1\}$ and

 $\{y\overrightarrow{C'}x_3x\overleftarrow{C'}z, y\overrightarrow{C'}x_3xx_4\overleftarrow{C'}z,
y\overrightarrow{C'}x_2xx_4\overleftarrow{C'}z\}$ with 

\hspace{2mm} $\{p_2+e_0+1+O_2+q_2,p_2+e_0+2+q_2,p_2+2+q_2\}$.

 \vskip 1mm

If any two of these three sets have different  path length sets, then we have at least $4=k+1$ $(y,z)$-paths of different lengths.
Otherwise,   we can check that $p_1=p_2=p$, $q_1=q_2=q$ and $O_1=O_2=O=e_0+1$. And we have three $(y,z)$-path lengths $p+2O+q,p+O+1+q,p+2+q$.
Since $|y\overrightarrow{C'}z|=p+O-1+q$, either we have $4=k+1$ $(y,z)$-paths of different lengths, or $O=3$ and $e_0=2$. In the second case, we have three $(y,z)$-path lengths $p+q+6,p+q+4,p+q+2$.

Since $xx_0x'\ss C'$,  $x'\notin \{y,z\}$, and either $|x\overrightarrow{C'}x'|=2$ or $|x\overleftarrow{C'}x'|=2$. 
 Using  the same arguments for $x'$, we only need to consider $x'$   satisfying  $s=1$ and $t=k=3$.    Suppose $x_1',x_2',x_3',x_4'$ follow each other on $C'$ in the forward direction.  If $|x\overrightarrow{C'}x'|=2$, then $|x\overleftarrow{C'}x_4'|=1$, and 
 we can find a $(y,z)$-path $y\overleftarrow{C'}x'x_4'\overleftarrow{C'}x_4xx_2\overrightarrow{C'}z$ of length $p+q+8$. If $|x\overleftarrow{C'}x'|=2$, then $|x\overrightarrow{C'}x_1'|=1$, and
 we can find a $(y,z)$-path $z\overrightarrow{C'}x'x_1'\overrightarrow{C'}x_1xx_3\overleftarrow{C'}y$ of length $p+q+8$. In both cases, we get $4=k+1$ $(y,z)$-paths of different lengths. 
 \bqed

\begin{lem}\label{C2k-1}
  Let $H$ be a graph obtained by adding  $2k-1$ chords $xx_1,\ldots,xx_{2k-1}$ in  an  even cycle $C_0$, where $x_1,\ldots, x_{2k-1}$ follow each other on $C_0$ in the forward direction. Let $C_{x,j}=x\overrightarrow{C_0}x_jx$ and $C_{x,j}'=x\overleftarrow{C_0}x_jx$. Then  $|L_e(H)|\geq k$. Moreover,  if $|L_e(H)|=k$ and  $H$ has $k$ consecutive even cycle lengths,  then
  
(1) we have $k$ even cycles $C_{x,j_1},\ldots,C_{x,j_{k-1}},C_0$ of lengths $4,\ldots,2k,2k+2$ for some $j_1<\cdots <j_{k-1}$ and $|C_{x,j_i}|=|C_{x,j_{k-i}}'|$;

(2)  $|x\overrightarrow{C_0}x_{j_1}|=|x\overleftarrow{C_0}x_{j_{k-1}}|=3$ and $|x_{j_i}\overrightarrow{C_0}x_{j_{i+1}}|=2$;

(3) $|x\overrightarrow{C_0}x_1|=|x\overleftarrow{C_0}x_{2k-1}|=2$.
\end{lem}

\noindent{\bf{Proof.}}
Since $C_0$ is even, 
$C_{x,j}$ and $C_{x,j}'$  are both even or   odd cycles. We say $x_j$ is even if both $C_{x,j}$ and $C_{x,j}'$ are even. Otherwise, $x_j$ is odd. Note that all $C_{x,j}$ and $C_{x,j}'$ are smaller than $C_0$. And we have at least either  $k-1$ even $x_j$ or $k$ odd $x_j$. If we    
have  $k-1$ even $x_j$, then $|L_e(H)|\geq k$. If  
there are $k$ odd $x_j$, denote them as $x_{j_1'},\ldots,x_{j_{k}'}$. Then we have an even cycle $xx_{j_1'}\overrightarrow{C_0}x_{j_i'}x$ for every $2\leq i\leq k$. Together with $C_0$,  $|L_e(H)|\geq k$. 

If $|L_e(H)|= k$, then we have 
 exactly $k-1$ even $x_j$ and $k$ odd $x_j$. Let $x_{j_1},\ldots,x_{j_{k-1}}$ be even. 
Then we have $k-1$ even cycles $C_{x,j_1},\ldots,C_{x,j_{k-1}}$ and $k-1$ even cycles $C_{x,j_1}',\ldots,C_{x,j_{k-1}}'$,  all of which are smaller than $C_0$. Since $|L_e(H)|=k$, $|C_{x,j_i}|=|C_{x,j_{k-i}}'|$ for every $1\leq i\leq k-1$. Then $|x\overrightarrow{C_0}x_{j_1}|=|x\overleftarrow{C_0}x_{j_{k-1}}|=O_1$ and $O_1$ is odd. 
If $H$ has $k$ consecutive even cycle lengths, then $|C_0|=|C_{x,j_{k-1}}|+2$, implying that $O_1=3$. Then $|C_{x,j_1}|=|C'_{x,j_{k-1}}|=4$.
And we see (1) and (2) hold.
Since $d(x)=2k+1$ and $|C_0|=2k+2$, we see (3) holds. 
\bqed

\begin{lem}\label{C2k-1odd}
Let $H$ be a graph obtained by adding  $2k-1$ chords $xx_1,\ldots,xx_{2k-1}$ in  an  odd cycle $C_0$, where $x_1,\ldots, x_{2k-1}$ follow each other on $C_0$ in the forward direction. Let $C_{x,j}=x\overrightarrow{C_0}x_{j}x$ and $C_{x,j}'=x\overleftarrow{C_0}x_{j}x$. Then $|L_e(H)|\geq k$ and $C_0$ is larger than  the longest even cycle. Moreover, if $|L_e(H)|=k$ and $H$ has $k$ consecutive even cycle lengths,   then    we have at least $k-1$ even cycles $C_{x,j}$ and $k-1$ even cycles $C'_{x,j}$, and  (1) or (2) holds:

(1) We have $k$ even cycles $C_{x,j_1},\ldots,C_{x,j_k}$ of lengths $4,\ldots,2k+2$ for some $j_1<\cdots <j_{k}$,  $|x\overrightarrow{C_0}x_{j_1}|=3$ and $|x_{j_i}\overrightarrow{C_0}x_{j_{i+1}}|=2$.

(2) We have $k$ even cycles $C_{x,j_1}',\ldots,C_{x,j_k}'$ of lengths $4,\ldots,2k+2$ for some $j_1>\cdots >j_{k}$, 
$|x\overleftarrow{C_0}x_{j_1}|=3$ and  $|x_{j_i}\overleftarrow{C_0}x_{j_{i+1}}|=2$.
\end{lem}

\noindent{\bf{Proof.}} Since $C_0$ is odd, one of $\{C_{x,j},C_{x,j}'\}$ is even, and the other is odd. Then we have at least either $k$ even $C_{x,j}$ or $k$ even $C_{x,j}'$, implying that $|L_e(H)|\geq k$ and $C_0$ is larger than the longest even cycle. If $|L_e(H)|= k$, then we have exactly $k$ even 
$C_{x,j}$ and  $k-1$ even 
$C_{x,j}'$, or $k-1$ even 
$C_{x,j}$ and  $k$ even 
$C_{x,j}'$. 
If we have $k$ even 
$C_{x,j}$,
let $j_1<\cdots <j_k$ be the indices such that  $C_{x,j_i}$ is even for every $1\leq i\leq k$. Let $|x\overrightarrow{C_0}x_{j_1}|=O_1$.
Since $xx_{j_1}\overrightarrow{C_0}x_{j_{i+1}}x$ is also an even cycle of length $|C_{x,j_{i+1}}|-O_1+1$ for every $1\leq i\leq k-1$, we get $|C_{x,j_i}|=|xx_{j_1}\overrightarrow{C_0}x_{j_{i+1}}x|=|C_{x,j_{i+1}}|-O_1+1$. Then $|x_{j_i}\overrightarrow{C_0}x_{j_{i+1}}|=O_1-1$ and 
 we have $k$ even cycle lengths $O_1+1,\ldots, k(O_1-1)+2$. If these even cycles have consecutive even lengths, we see (1) holds.
If we have $k$ even $C_{x,j}'$, we can deduce (2) by the same arguments. \bqed

\vskip 5mm
\noindent{\bf\large Acknowledgements}
\vskip 2mm

The author sincerely thanks Yulong Niu  at Institut de la Vision, Sorbonne Université for helpful discussions, and  Anders Yeo at South Denmark University for carefully reading and providing  helpful comments. 

Wang was supported by National Key R\&D Program of China under grant number 2023YFA1010202 and NSFC under grant number 12401447.

\end{document}